\documentclass[final,3p,11pt,A4]{elsarticle}

\usepackage{graphicx}

 \usepackage{subfig}
 \usepackage[utf8]{inputenc}
 \usepackage{amsmath,amssymb,amsfonts,nomencl,mathtools} 
 \usepackage{booktabs}
 \usepackage{algorithmicx}

 \usepackage{amstext,amsmath,amssymb}
 \usepackage{algorithm}
 \usepackage{hyperref}
 \usepackage{algpseudocode}
 \usepackage{xcolor}
 \usepackage{natbib}
 
 \usepackage{multirow}
 
 \usepackage{tikz}
\usepackage{pgfplots}

\usepackage{appendix}

\usepackage{natbib}
\biboptions{authoryear}

\newcommand{\NULL}{\textsc{null}}
\newcommand{\Tau}{\mathcal{T}}

\newcommand{\pbs}{\pi^{bs}}

\DeclareMathOperator*{\argmax}{arg\,max}

\usepackage[final,markup=none,commentmarkup=footnote]{changes}
\usepackage{color}
\setcommentmarkup{\footnote{{\bf #2 \arabic{authorcommentcount}:} \ifthenelse{\isColored}{\color{authorcolor}}{}#1}}
\definechangesauthor[color=orange]{YS}
\definechangesauthor[color=maroon]{AE}
\definechangesauthor[color=blue]{XL}

\definechangesauthor[color=red]{DT}

\definecolor{marine}{RGB}{0,32,96}
\definecolor{navy}{RGB}{0,0,128}
\definecolor{maroon}{RGB}{128,0,0}
\definecolor{olivegreen}{RGB}{85,107,47}
\definecolor{gray}{RGB}{102,102,102}
\definecolor{green}{RGB}{131,198,210}
\definecolor{blue}{rgb}{0, 0.4470, 0.7410}
\definecolor{skyblue}{rgb}{0.3010, 0.7450, 0.9330}

\usepackage{xcolor,colortbl}
\definecolor{Gray}{gray}{0.90}
\definecolor{LightCyan}{rgb}{0.88,1,1}
\definecolor{wheat}{rgb}{0.96, 0.87, 0.7}
\definecolor{teagreen}{rgb}{0.82, 0.94, 0.75}

\newcolumntype{g}{>{\columncolor{Gray}}l}
\newcolumntype{o}{>{\columncolor{orange}}l}
\newcolumntype{s}{>{\columncolor{LightCyan}}l}
\newcolumntype{w}{>{\columncolor{wheat}}l}
\newcolumntype{t}{>{\columncolor{teagreen}}l}

\pgfplotsset{every tick label/.append style={font=\tiny}}
\pgfplotsset{
	box plot width/.initial=4em,
	box plot/.style={
		/pgfplots/.cd,
		black,
		only marks,
		mark=-,
		mark size=\pgfkeysvalueof{/pgfplots/box plot width},
		/pgfplots/error bars/.cd,
		y dir=plus,
		y explicit,
	},
	box plot box/.style={
		/pgfplots/error bars/draw error bar/.code 2 args={%
			\draw [line width=0.20mm]  ##1 -- ++(\pgfkeysvalueof{/pgfplots/box plot width},0pt) |- ##2 -- ++(-\pgfkeysvalueof{/pgfplots/box plot width},0pt) |- ##1 -- cycle;
		},
		/pgfplots/table/.cd,
		y index=2,
		y error expr={\thisrowno{3}-\thisrowno{2}},
		/pgfplots/box plot
	},
	box plot top whisker/.style={
		/pgfplots/error bars/draw error bar/.code 2 args={%
			\pgfkeysgetvalue{/pgfplots/error bars/error mark}%
			{\pgfplotserrorbarsmark}%
			\pgfkeysgetvalue{/pgfplots/error bars/error mark options}%
			{\pgfplotserrorbarsmarkopts}%
			\path ##1 -- ##2;
		},
		/pgfplots/table/.cd,
		y index=4,
		y error expr={\thisrowno{2}-\thisrowno{4}},
		/pgfplots/box plot
	},
	box plot bottom whisker/.style={
		/pgfplots/error bars/draw error bar/.code 2 args={%
			\pgfkeysgetvalue{/pgfplots/error bars/error mark}%
			{\pgfplotserrorbarsmark}%
			\pgfkeysgetvalue{/pgfplots/error bars/error mark options}%
			{\pgfplotserrorbarsmarkopts}%
			\path ##1 -- ##2;
		},
		/pgfplots/table/.cd,
		y index=5,
		y error expr={\thisrowno{3}-\thisrowno{5}},
		/pgfplots/box plot
	},
	box plot median/.style={
		/pgfplots/box plot
	}
}
\pgfplotsset{yticklabel style={text width=1.0em,align=right}}
\pgfplotsset{compat=1.12}

\usepackage{lipsum}
\makeatletter
\def\ps@pprintTitle{%
	\let\@oddhead\@empty
	\let\@evenhead\@empty
	\def\@oddfoot{}%
	\let\@evenfoot\@oddfoot}
\makeatother

\begin{document}

\begin{frontmatter}
	
	\title{Instance Space Analysis for the Car Sequencing Problem}
	\author[label1]{Yuan Sun}
	\ead{yuan.sun@monash.edu}
	\author[label1]{Samuel Esler}
	\ead{sam.p.esler@gmail.com}
	\author[label2]{Dhananjay Thiruvady}
	\ead{dhananjay.thiruvady@deakin.edu.au}	
	\author[label1]{Andreas Ernst}
	\ead{andreas.ernst@monash.edu}
	\author[label3]{Xiaodong Li}
	\ead{xiaodong.li@rmit.edu.au}
	\author[label2]{Kerri Morgan}
	\ead{kerri.morgan@deakin.edu.au}	
	
	\address[label1]{School of Mathematical Sciences, Monash University, Clayton, 3800, Victoria, Australia}
	\address[label2]{School of Information Technology, Deakin University, Geelong, 3217, Victoria, Australia}
	\address[label3]{School of Computing Technologies, RMIT University, Melbourne, 3001, Victoria, Australia}
		
	
	
	\begin{abstract}
		We investigate an important research question for solving the car sequencing problem, that is, which characteristics make an instance hard to solve? To do so, we carry out an instance space analysis for the car sequencing problem, by extracting a vector of problem features to characterize an instance. \added[]{In order to visualize the instance space, the feature vectors are projected onto a two-dimensional space using dimensionality reduction techniques.} The resulting two-dimensional visualizations provide new insights into the characteristics of the instances used for testing and how these characteristics influence the behaviours of an optimization algorithm. This analysis guides us in constructing a new set of benchmark instances with a range of instance properties. We demonstrate that these new instances are more diverse than the previous benchmarks, including some instances that are significantly more difficult to solve. \added[]{We introduce two new algorithms for solving the car sequencing problem and compare them with four existing methods from the literature.} \added[]{Our new algorithms are shown to perform competitively for this problem but no single algorithm can outperform all others over all instances.} This observation motivates us to build an algorithm selection model based on machine learning, to identify the niche in the instance space that an algorithm is expected to perform well on. Our analysis helps to understand problem hardness and select an appropriate algorithm for solving a given car sequencing problem instance. 		
		
	\end{abstract}
	
	
\end{frontmatter}

\section{Introduction}
\label{intro}
The production of cars in Europe rely on a ``built-to-order'' system, where cars with specific options need to be assembled in relatively short periods of time. The cars to be assembled are typically diverse in respect to their required options, and hence, the assembly process is necessarily slow. Since 1993, the car manufacturer, Renault, has been tackling a variant of ``car sequencing'' with the aim of assembling a range of diverse cars as efficiently as possible \citep{Liris-2656}. In 2005, the subject of the ROADEF challenge was the same car sequencing problem \citep{nguyen_roadef}, which is one of finding an optimal order, or sequence, in which cars are to be processed on a production line.  Each car has a set of options such as paint colour, air conditioning, etc. Each option has an associated constraint that gives the maximum number of cars that can have this option applied in any contiguous subsequence. The car sequencing problem is known to be NP-hard \citep{KIS2004331}.

\added[]{There are many variants of car sequencing problems that all share common characteristics but with small differences in the objective function or constraint. A comprehensive study of all variants is beyond the scope of a single paper. Here we will focus on one} popular optimization variant of the car sequencing problem which minimizes the violation of constraints on the occurrence of options within subsequences \citep{bautista08, Thiruvady11, Thiruvady2014a}. The first study in this direction was conducted by \citet{bautista08}, who showed that a beam search approach can be very effective. Following this, \citet{Thiruvady11} showed that a hybrid  of beam search, ant colony optimization (ACO) and constraint programming can substantially improve solution qualities while reducing run-times. \citet{Thiruvady2014a} proposed a hybrid of Lagrangian relaxation and ACO, and they showed the hybrid can find excellent solutions and at the same time very good lower bounds. A recent study in this direction is by \citet{thiruvady2020large}, who showed that a large neighbourhood search (LNS) finds high quality solutions in relatively low run-times. \added[]{Other methods used for this type of problem include beam search \citep{golle2015iterative} and a hybrid neighbourhood search method combining tabu search and LNS~\citep{zhang2018hybrid}.}


In the literature, there are three commonly used benchmarks \citep{gent98, gravel2005review,perron2004combining}. The benchmark instances from \citet{perron2004combining} consisted of up to 500 cars with relatively high utilisation of options, and were considered as the ``hardest'' set of problem instances available to test algorithms against. Though, since then, studies \citep{Thiruvady11,Thiruvady2014a,thiruvady2020large} have demonstrated different solvers (including those based on mixed integer programming (MIP), constraint programming and metaheuristics) have variable performance, and hence it is unclear as to which set of problem instances are in fact the most complex. Moreover, the specific characteristics of hard problem instances have not clearly been identified. 


\added[]{The aims of this study are twofold: (a) to address the gap of a lack of knowledge as to which characteristics of the car sequencing problem make finding high quality solutions hard, and (b) to use machine learning to build an algorithm selection model.} The first aim is achieved by performing an instance space analysis for the problem. The methodology of instance space analysis was originally proposed by \citet{smith2014towards}, and has been successfully applied to many problems including combinatorial optimization \citep{smith2014towards}, continuous optimization \citep{munoz2017performance}, machine learning classification \citep{munoz2018instance}, regression \citep{munoz2021instance} and facial age estimation \citep{smith2020revisiting}. However, it still requires problem-specific design to apply this general approach to the specific problem of car sequencing, such as feature extraction and instance generation. 


\added[]{To carry out an instance space analysis,} we extract a vector of problem features to characterise a car sequence problem instance, and use principal component analysis  to reduce the dimension of feature vectors into two. By doing so, an instance is projected to be a point in the two-dimensional instance space, allowing a clear visualization of the distribution of instances, their feature values and algorithms' performance. Our analysis reveals two problem characteristics that make a car sequencing instance hard to solve, i.e., high utilisation of options and large average number of options per car class. Further, this analysis allows us to develop a problem instance generator, which can generate problem instances that can be considered extremely complex irrespective of the type of solver. 

\added[]{The second aim, building an algorithm selection model, is achieved by using machine learning} \citep{rice1976algorithm,smith2009cross,munoz2015algorithm} to identify which algorithm is most effective in a particular region of the instance space. We select three state-of-the-art algorithms from the literature as our candidates: the Lagrangian relaxation and ACO hybrid \citep{Thiruvady2014a} and two LNS algorithms proposed by \citet{thiruvady2020large}. Additionally, given the vast improvement in commercial solvers, we investigate three algorithms not previously successful or attempted in the literature. These algorithms are MIP, MIP based on lazy constraints, and an adaptive LNS that builds on algorithms studied in \citet{thiruvady2020large}.




\noindent 
In summary, the contributions of this study are: 
\begin{enumerate}
    \item \added[]{Introducing two new optimization algorithms for the car sequencing problem: an adaptive LNS method and a MIP method with lazy constraints. We also demonstrate for the first time that with suitable tuning of parameters, modern MIP solvers are competitive for this problem.}
    \item Performing an instance space analysis for the car sequence problem and successfully identifying the characteristics that make an instance hard to solve.  
    \item Designing a benchmark to allow systematical generation of problem instances with controllable properties and showing that they complement the existing benchmark instances well in the instance space.
    \item Developing an algorithm selection model based on machine learning, which shows that our adaptive LNS and the MIP solver based on Gurobi become the new state-of-the-art for solving the car sequencing problem. Further, the adaptive LNS performs the best for solving hard problem instances, while the MIP solver is the best for medium-hard instances. 
\end{enumerate}
\added[]{This study is different to many others published in the literature in not simply identifying one new or improved algorithm but advancing the state-of-the-art while giving insight into what types of algorithms are best suited for different types of instances. The analysis of the instance space and new benchmark instances are expected to pave the way for future research in designing better algorithms for car sequencing problems.}


The remainder of this paper is organized as follows. In Section~\ref{sec::instance space analysis}, we introduce the methodology of algorithm selection and instance space analysis. In Section~\ref{sec:prob}, we describe the car sequencing problem and present the MIP formulation used. \added[]{Section~\ref{sec:algorithms} introduces the solution algorithms for the car sequencing problem, and Section~\ref{sec::features and instances} describes feature extraction and instance generation.} Section~\ref{sec::experiments} presents and analyses the experimental results, and the last section concludes the paper.

\section{Algorithm Selection and Instance Space Analysis}
\label{sec::instance space analysis}


The algorithm selection framework was originally proposed by \citep{rice1976algorithm}. It has four key components:
\begin{itemize}
    \item Problem space ($P$) consists of instances from the problem of interest.
    \item Algorithm space ($A$) contains algorithms available to solve the problem.
    \item Feature space ($F$) contains measurable features extracted to characterise a problem instance.
    \item Performance space ($Y$) consists of measures to evaluate the performance of algorithms when solving the instances, e.g., the optimality gap generated by an algorithm within a time limit.   
\end{itemize}
For an instance $x\in P$, a set of features $f(x)$ can be first extracted. The goal of algorithm selection is to learn a mapping from the feature vector $f(x) \in F$ to an algorithm $\alpha^* \in A$, such that the algorithm $\alpha^*$ performs the best on the instance $x$ when compared to other algorithms in $A$. This is typically achieved by training a machine learning model on the data composed of the feature vectors of training instances and the performances of candidate algorithms on the training instances. \added[]{Algorithm selection is closely related to meta-learning \citep{vanschoren2018meta} and auto machine learning \citep{he2021automl}, which aims to select the best parameter settings for machine learning models.}

\added[]{Algorithm selection is often beneficial for problem solving. For example the algorithm selection model, SATzilla-07, was able to solve more than 20\% instances in the 2007 SAT competition than any of its component solver \citep{xu2008satzilla}. Apart from the SAT problem, algorithm selection has been demonstrated successful on quantified Boolean formulas \citep{pulina2009self}, constraint solving \citep{o2008using}, answer set programming \citep{hoos2014claspfolio}, continuous black-box optimization \citep{kerschke2019automated} and clustering \citep{wang2020efficiency}, to name a few. Some of those algorithm selection models were collected in the ASlib benchmark library \citep{bischl2016aslib}, and more comprehensive literature reviews are available in \citep{smith2009cross,munoz2015algorithm,kerschke2019automated}.}


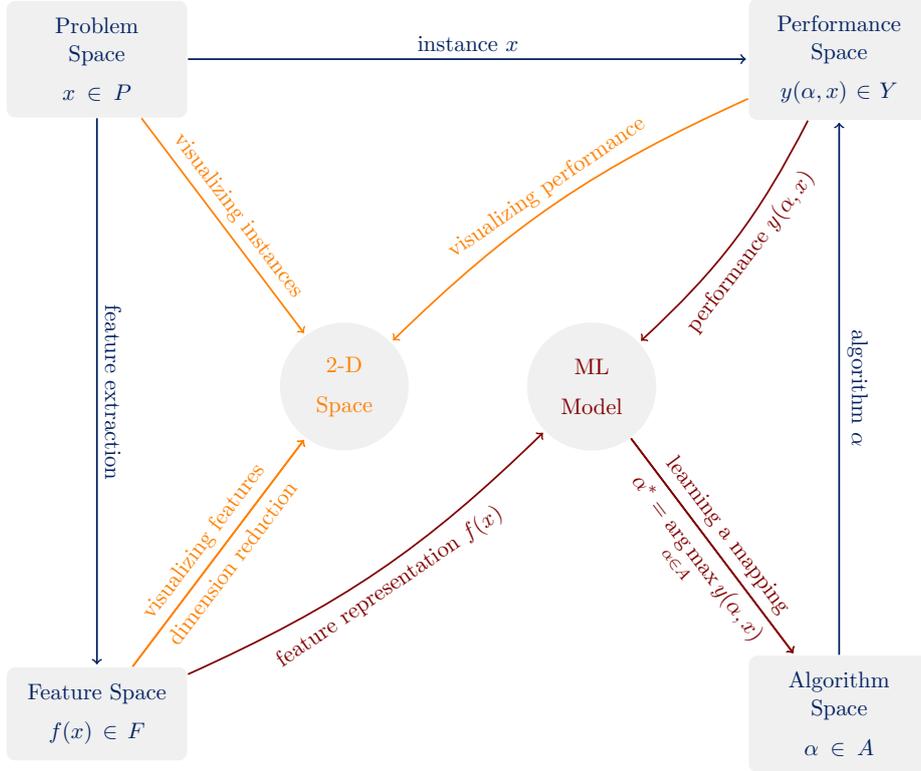
\begin{figure}[t]
	\centering
	\resizebox{0.75\textwidth}{!}{
	\begin{tikzpicture}[shorten >=1pt,auto,node distance=1cm, thick, main node/.style={color = marine, rectangle, fill=gray!10, rounded corners, inner sep=8pt, text width=2.5cm, align=center}, node1/.style={color = orange, circle, fill=gray!10, inner sep=8pt, align=center}, node2/.style={color = maroon, circle, fill=gray!10, inner sep=8pt, align=center}, scale = 1.4]
	\node[main node] (1) at (0,0) {Feature Space \\[0.5em] $f(x) \in F$};
	\node[main node] (2) at (9,0) {Algorithm Space \\[0.5em] $\alpha \in A$};
	\node[main node] (3) at (0,8) {Problem Space \\[0.5em] $x \in P$};
	\node[main node] (4) at (9,8) {Performance Space \\[0.5em] $y(\alpha, x) \in Y$ };
	\node[node1] (5) at (3,4) {2-D \\[0.5em] Space};
	\node[node2] (6) at (6,4) {ML \\[0.5em]  Model};

	\path[every node/.style={color = marine}]
	(3) [color = marine, line width = 0.30mm, ->, sloped] edge node [] {instance $x$} (4)
	(2) [color = marine, line width = 0.30mm, ->, sloped] edge node [] {algorithm $\alpha$} (4)
	(3) [color = marine, line width = 0.30mm, ->, sloped] edge node [] {feature extraction} (1);
	
	\path[every node/.style={color = marine}]
	(1) [color = orange, line width = 0.30mm, ->] edge node [sloped, below, color = orange] {dimension reduction} (5)
	(1) [color = orange, line width = 0.30mm, ->] edge node [color = orange, sloped, above] {visualizing features} (5)
	(3) [color = orange, line width = 0.30mm, ->] edge node [color = orange, sloped, above] {visualizing instances} (5)
	(4) [color = orange, line width = 0.30mm, ->, bend right = 10] edge node [color = orange, sloped, above] {visualizing performance} (5);
	
	\path[every node/.style={color = marine}]
	(1) [color = maroon, line width = 0.30mm, ->, bend right = 10] edge node [sloped, below, color = maroon] {feature representation $f(x)$} (6)
	(4) [color = maroon, line width = 0.30mm, ->, bend left = 10] edge node [sloped, below, color = maroon] {performance $y(\alpha, x)$} (6);
	\path[every node/.style={color = marine}]
	(6) [color = maroon, line width = 0.30mm, ->] edge node [sloped, above, color = maroon] {learning a mapping} (2)
	(6) [color = maroon, line width = 0.30mm, ->] edge node [sloped, below, color = maroon] {$\alpha^* = \argmax\limits_{\alpha \in A} y(\alpha, x)$} (2);
	\end{tikzpicture}
	}
	\caption{The methodological framework of algorithm selection and instance space analysis. \added[]{Note that the reduced 2-D space is only used for visualization, and the algorithm selection models are built on the feature space $F$.}}
	\label{methodology}
\end{figure}

\citet{smith2014towards} further extended Rice's algorithm selection framework by including a visualization component (See Fig.~\ref{methodology}). A significant advance of the extended algorithm selection framework is that it provides great insights into the strengths and weaknesses of an algorithm as well as problem hardness. A two-dimensional instance space can be created and visualised by using dimension reduction techniques such as principal component analysis to reduce the dimension of feature vectors to two. Importantly, the performance of an algorithm can also be visualized in the two-dimensional space, clearly showing in which region of the instance space the algorithm performs well or badly. Finally, by visualizing the values of features in the same space, one can identify the features that are mostly correlated with the performance of algorithms, indicating which features make an instance hard to solve. 

\added[]{Despite the existence of many studies in this area, the existing algorithm selection models are all problem-specific, and thus cannot be directly applied to car sequencing problems. In the following, we will perform an instance space analysis and build algorithm selection models for the car sequencing problem.}

\section{The Car Sequencing Problem} \label{sec:prob}
In this section, we describe the optimization version of the car sequencing problem, originally introduced by \citet{bautista08}. 

Consider $D$ cars, each of which belongs to a class $C=\lbrace c_{1}, \ldots, c_{K}\rbrace$, where $K$ is the number of classes. The cars require one or more of a set of options $O=\lbrace o_{1}, o_{2}, \ldots, o_{O}\rbrace$ that need to be installed. All cars within a class must be fitted with the same options. For this purpose,  each car $1,\ldots,D$ has an indicator vector associated with it  
\begin{align}\label{eqn:rdef}
 \vec{r_{i}}[j]=\begin{cases}
               1&\text{if car } i \text{ requires option }j,\\
               0&\text{otherwise}
              \end{cases}\qquad \text{ for }1\leq j \leq O.
\end{align}
Cars in the same class have the same indicator vector, or alternatively two cars $i$ and $k$ are in the same class if $\vec{r_{i}}=\vec{r}_{k}$. Moreover, the number of cars in class $j$ is denoted as $d_{j}$.  


Typically, each station that installs the options along the assembly line has a certain capacity, which if not met, will slow down the assembly line. Hence, to satisfy aspirational targets, we introduce for each option $p_{j}$ and $q_{j}$ ($j=1,\ldots,O$), where in any subsequence of $q_j$ cars, there should only be $p_j$ cars that require option $j$. 

Any sequence or ordering of the cars is a potential solution to the problem. We denote such a sequence by $\pi$, and $\pi_{k}$ is used to represent car class $c$ in position $k$ of $\pi$. 
As any permutation between cars in the same class produces an equivalent solution, we eliminate such symmetries in the problem in order to reduce the search space.  

A sequence that satisfies all subsequence constraints for all options can be considered as ideal. Though, such a sequence can be hard to find, and hence the optimization version of the problem focuses on modulating option utilisation across the sequence. The measures that capture the modulation of options are presented next.

\subsection{Assignment measures} \label{sec:ass_measures}
An option $o_{i}$ with subsequence length $q_{i}$ will have $D-q_{i}+1$ subsequences. Those subsequences which have used $o_{i}$ in excess of the allowable $p_{i}$ options, can be computed as follows:  
\begin{align}
 y_{ij}(\pi)=\max \left \lbrace 0, - p_{i} + \sum_{k=u_{i}(j)}^{j} r_{\pi(k)}[i] \right \rbrace
\end{align}
where $u_{i}(j)=\max(1, j+1-q_{i})$. Here $u_{i}(j)$ is the starting position of a subsequence $q_{i}$, which ends at position $j$.  The inequality $j<q_{i}$ implies that the subsequence starts at position $1$, resulting in the subsequence size being less than $q_{i}$.

We now define \textit{upper over-assignment} ($uoa$) \citep{bautista08}, which computes the total number of times within a sequence $\pi$ that an option $o_{i}$ violated the limit $p_{i}$: 
\begin{align}
 uoa(\pi)=\sum_{i=1}^{O}\sum_{j=p_{i}}^{D}a_{ij}y_{ij}(\pi)\label{eqn1}
\end{align}
where $a_{ij}$ is a penalty imposed for exceeding $p_j$ in subsequence positions $j+1-q_{i}$ to $j$ (if $j<q_{i}$, positions $1$ to $j$). If $p_j$ is never violated, then $uoa(\pi)=0$, and $\pi$ is considered as a feasible solution.  

The \textit{upper under-assignment} ($uua$) is a second measure that determines the number of times, and by how much, the utilisation of options within the subsequences are under $p_j$:  
\begin{align}
 uua(\pi)&=\sum_{i=1}^{O}\sum_{j=p_{i}}^{D}b_{ij}z_{ij}(\pi)\label{eqn2}\\
 z_{ij}(\pi)&=\max \left \lbrace 0, p_{i} - \sum_{k=u_{i}(j)}^{j} r_{\pi(k)}[i] \right \rbrace,
\end{align}
where $b_{ij}$ is a penalty imposed for having fewer than $p_{i}$ cars in subsequence positions $j+1-q_{i}$ to $j$ (if $j<q_{i}$, positions $1$ to $j$). 

Penalties for over assignment are typically different to that of under assignment, and the penalties in both cases can vary across a sequence. If the penalties were set as constants, it would suffice to measure only over assignments (and hence over assignment penalties), as the total increase in over assignments leads to an equivalent decrease in under assignments in other parts of the sequence.

\subsection{Mixed Integer Programming Model}\label{sec:IP}

\citet{bautista08} provided an efficient MIP formulation for this problem. We use this formulation for the MIP-based methods that follow in this study. 

As introduced earlier, $y$ and $z$ are chosen as variables (Section~\ref{sec:ass_measures}). Moreover, we define binary variables $x$ where $x_{it} = 1$ if car class $i$ is chosen in position $t$. The objective is to minimize $uoa(\pi)+uua(\pi)$ corresponding to Equations \eqref{eqn1} and \eqref{eqn2}. The problem can be defined as follows:
\begin{alignat}{3}
 &&\text{Min.} \quad\sum_{i=1}^{O}\sum_{j=p_{i}}^{D}&\left(a_{ij}y_{ij}+b_{ij}z_{ij}\right)&&\label{goal}\\
 \noalign{\text{Subject to}} \nonumber \\
 &&\sum_{i=1}^{K}x_{it}&=1 && t \in [1,\ldots,D] \label{someClassAtPosition} \\ 
  &&\sum_{t=1}^{D}x_{it}&=d_{i}&&i \in [1,\ldots,K] \label{allCars}\\
  &&\sum_{i=1}^{K}\sum_{u_j(t)
  }^{t}\vec{r}_i[j] \cdot x_{ik}&=p_{j}+y_{jt}-z_{jt}&\qquad &\forall j\in[1,\ldots,O], \forall t \in [p_{j},\ldots,D]\label{satOptions}\\
&&z_{jt}, y_{jt}&\geq 0 && \forall j \in [1,\ldots,O], \forall t \in [p_{j},\ldots,D]\\
&&x_{it}&\in [0,1]&&\forall i \in [1,\ldots,K], \forall t\in[1,\ldots,D].
\end{alignat}

This formulation differs slightly from  the discussion of the problem in the previous section. Observing that two cars of the same class are interchangeable (i.e., swapping two cars of the same class leads to the same solution), this formulation focuses on selecting car classes in each position. It has the added advantage of reducing the solution space by removing symmetries. 


\section{Algorithms}\label{sec:algorithms}

\added[]{A large number of solution methods have been proposed to tackle at least some variants of the car sequencing problem, including exact methods such as constraint programming \citep{dincbas88}, MIP \citep{gravel2005review}, meta-heuristics \citep{gottlieb03,SUN201893} and hybrid methods \citep{kichane08}. A comprehensive review of this literature is beyond the scope of this paper. Here, we focus on techniques for solving the optimization variant of the car sequencing problem, including Lagrangian relaxation and ant colony optimization (LR-ACO) \citep{Thiruvady2014a}, and large neighbourhood search (LNS) \citep{thiruvady2020large}. Further, we introduce a new adpative LNS algorithm and investigate a MIP method with lazy constraints.}



\subsection{Lagrangian Relaxation and Ant Colony Optimization} \label{sec:lraco} 

An effective Lagrangian relaxation procedure from the MIP model (Section~\ref{sec:IP}) was proposed by \citet{Thiruvady2014a}. Among the different sets of constraints that can be relaxed, they showed that Constraint~\eqref{someClassAtPosition} leads to the most effective Lagrangian relaxation. Hence, the resulting problem consists of all other constraints and the objective is modified as follows:
\begin{equation}\label{one}
	LLR(\lambda) = \min \sum_{j=1}^O \sum_{t=1}^D (a_{jt}  y_{jt} + b_{jt}   z_{jt}) + \sum_{t=1}^D \lambda_t  \left(\sum_{i=1}^K x_{it}\,-1\right).
\end{equation}

Solving this relaxation leads to a solution that is potentially infeasible. Specifically, certain positions ($i$ in $x_{it}$) may have no car class assigned and other positions may have two or more cars assigned. To rectify this infeasibility, Lagrangian multipliers ($\lambda$) are introduced, which penalise positions in the sequence where either zero or more than one car have been assigned. The multipliers are updated via subgradient optimization \citep{bertsekas95}, which aims to find the ``ideal'' set of multipliers over a number of iterations. The ACO algorithm \citep{dorigo04} is incorporated within Lagrangian relaxation \added[]{as a repair method}, leading to the LR-ACO heuristic. The flowchart of the LR-ACO method is shown in Fig.~\ref{fig:lraco_flowchart} and its pseudocode and associated details are presented in Appendix~A. 

\begin{figure}[!t]
    \centering
    \includegraphics[scale=0.50]{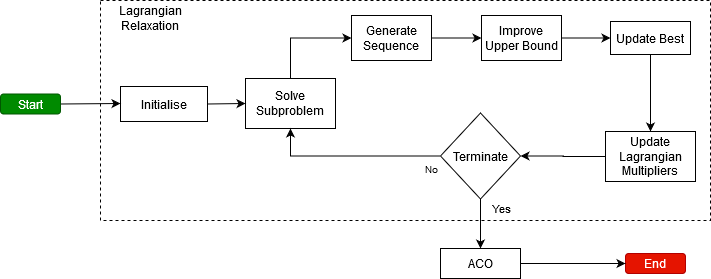}
    \caption{A high-level view of the LR-ACO heuristic, where Lagrangian relaxation is first applied, and 
    ACO is then used as a repair method.}
    \label{fig:lraco_flowchart}
\end{figure}

\subsection{Large Neighbourhood Search}\label{sec:lns}

The LNS of \citet{thiruvady2020large} is used as the basis for the LNS implementations in this study. At a high-level the procedure works by starting with a (good) feasible solution, considering subsequences within the sequence, solving a MIP of the subsequences, and repeating this process until some termination criterion is met. The motivation for this approach is that solving the original MIP can be very time consuming, and hence solving smaller part of the original problem and piecing them together can be achieved efficiently. 


\begin{figure}[!t]
    \centering
    \includegraphics[height=4cm,width=1.0\textwidth]{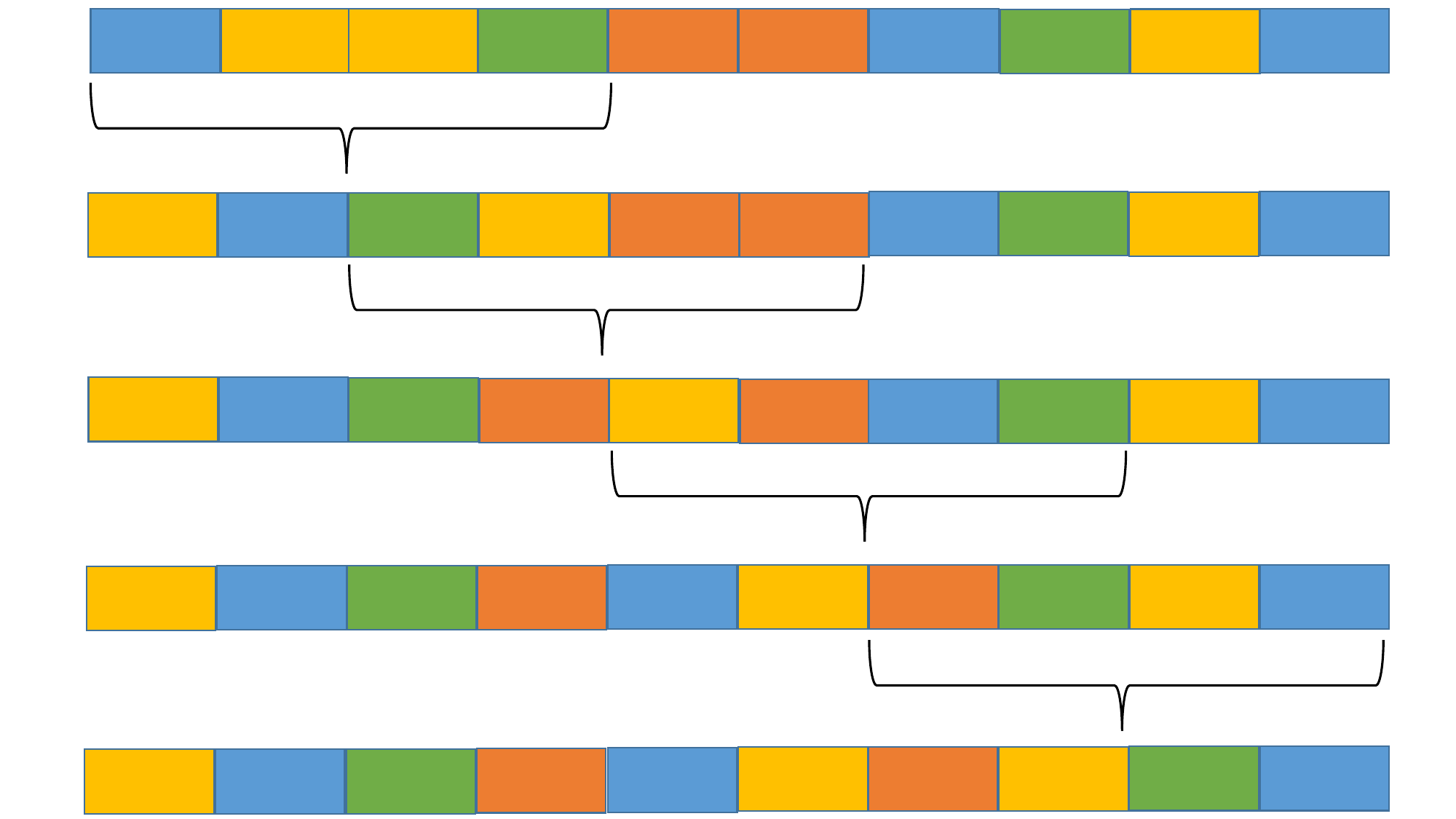}
    \caption{An illustration of one iteration of the LNS method. In this example, there are ten cars and four car classes, indicated by the colours: blue, green, yellow and orange. The first four cars are considered in the MIP, and solved leading to swapping all four car classes. Following this step, there is a shift of two cars, and the next four cars are considered to be solved. The search in this iteration continues until reaching the end of the sequence.}
    \label{fig:lns_example}
\end{figure}

The procedure works as follows. First, a high quality solution is obtained in a short time-frame. Following this step, subsequences within the solution are considered sequentially. For each subsequence, a MIP (Section~\ref{sec:IP}) is generated (i.e., subproblem) that fixes all variables in the original solution, but leaves those solution components associated with the subsequence free. The MIP is then solved, possibly to a time limit or gap for efficiency purposes.\footnote{The study by ~\cite{thiruvady2020large} used a gap of 0.5\% obtained via preliminary test.} If the subsequence sizes are small, all subproblems can be solved efficiently. In some cases large subproblems can also be efficiently solved. In one iteration of the algorithm, the search progresses by solving a subsequence, defining a new subsequence by moving forward a few positions, and continuing this process until the end of the subsequence is reached (see Fig.~\ref{fig:lns_example} for an example). Key to the success of the method is to move forward a few positions by ensuring at least some overlap with the previous subsequence. The above process is repeated for multiple iterations until the termination criterion is reached. The pseudocode of LNS and its associated details are presented in Appendix B. 

The size of the subsequence (i.e., window size) has a large impact on the performance of LNS. In our experiments, we will test two versions of LNS with different window sizes $w$: (1) 10-LNS with $w=10$, and (2) LCM-LNS with $w$ equal to the lowest common multiple of $q$ values. 

\subsection{A \added[]{New} Adaptive LNS Approach}\label{sec:improved_lns}

We make three modifications to the LNS algorithm proposed by \citet{thiruvady2020large} to improve its performance\added[]{, creating a new adaptive LNS}. 

Firstly, we note that it is a waste of computational resources when solving consecutive subsequences of cars repeatedly. Considering a car sequence $\pi$ to be improved, if a number of initial subsequences of $\pi$ are already solved to optimality, there is no need to re-optimize those subsequences in the following iterations. To address this issue, we propose to optimize \emph{non-consecutive} subsequences of $\pi$ so that different subsequences are optimized in each iteration of the algorithm. More specifically, we generate a random permutation ($P$) of integers from $1$ to $D$ in each iteration, and solve non-consecutive subsequences of $\pi$ indexed from $P[\hat{s}]$ to $P[\hat{s}+w]$, where $w$ is the window size, and $\hat{s}$ is the starting point, which increments by an interval of $w/2$ to generate multiple subsequences. Note that the randomization of $P$ has the effect that any large neighbourhood consisting of $w$ cars may be searched, rather than only subsequences of $w$ consecutive cars as was the case in the original LNS.

Secondly, the original LNS algorithm uses a fixed window size $w$ for an instance, i.e., setting $w$ to a constant (e.g., ten) or the lowest common multiple of $q$ values. However, we observe that for some instances, the algorithm can solve subsequences of small sizes very quickly and may get stuck in a local optimum in an early stage of the search. To address this potential issue, we introduce an adaptive setting for the window size, that is, starting with a window size of thirty and increasing the window size by one if the improvement of the optimality gap in an iteration is less than 0.005. By doing so, the algorithm can escape from a local optimum by using a larger window size. If computational time allowance is sufficient and the window size is eventually increased to $D$, solving the subsequence is equivalent to solving the whole sequence (i.e., the original problem).

In addition, given the vast improvement in commercial solvers, we use Gurobi to solve the MIP within a cutoff time to generate an initial feasible solution for the LNS algorithm, instead of using LR-ACO. This modification makes the algorithm simpler and more straightforward to implement.

\subsection{MIP with Lazy Constraints} \label{sec:lazy}
In similar fashion to Lagrangian relaxation, the lazy constraint method relaxes one or more of the sets of constraints. As we have previously motivated, solving the relaxed problem is typically time efficient, and the aim in the lazy method is to exploit this potential advantage. The key idea behind this procedure is that solutions typically violate a small number of constraints within the constraint set. If chosen carefully, it is possible to find high quality solutions with little overhead.  

The procedure works as follows. Prior to the execution of the algorithm, a constraint in the MIP formulation is chosen to be relaxed. The algorithm initializes the MIP model and executes without considering the chosen constraint. When a feasible solution is found, it is tested for feasibility against the constraint set that is relaxed. If violated, the specific constraint is introduced to the model as a ``lazy'' constraint (note, not the whole constraint set), and the procedure executes again. This process repeats until a provably optimal solution (to the original problem) is found or until a time limit expires. 

For the purposes of this study, and in preliminary experimentation, we found that Constraints \eqref{satOptions} are the most effective as implemented within the lazy method. This finding is not surprising, since relaxing this constraint set leads to a MIP that is equivalent to an assignment or bi-partite matching problem \citep{Manlove2013}. These problems are known to be pseudo-polynomial, and can hence be solved efficiently. Note, there is a potential disadvantage to this approach, which is that the link between the $x$ and $y$/$z$ variables are broken. It means for specific problem instances, an initial feasible solution can be hard to find (i.e. too many lazy constraints are needed to find a feasible solution).

\section{Feature Extraction and Instance Generation}\label{sec::features and instances}

\subsection{Features}
\label{features}

\added[]{We seek to determine the features of the car sequencing problem that will allow generating hard problem instances, i.e., ones that are considered hard for existing and newly proposed methods. There are certain features that are shown to lead to problem hardness as identified in previous studies (e.g. number of cars), and we also consider other features which could lead to creating hard instances \citep{gent1999csplib,gravel04,perron04}. The following are the features chosen using past studies as a guide. 
}
\begin{itemize}
  \item[] 1. Number of cars
  \item[] 2. Number of options
  \item[] 3. Number of car classes
  \item[] 4-7. Option utilisation (minimum, average, maximum, standard deviation)
  \item[] 8. Average number of options per car class
  \item[] 9-12. $p/q$ ratio (minimum, average, maximum, standard deviation)
  \item[] 13. Lowest common multiple of $q$ values
  \item[] 14. Standard deviation of car class populations
\end{itemize}

Note, the utilisation of an option $o$ is defined as: $\mu_o(\sum_{i = 1}^{D}r_i)/D$, where $r_i = 1$ if and only if car $i$ requires option $o$ and $\sum_{i = 1}^{D}r_i$ is the total number of cars that require option $o$, and $\mu_o(t)$ is the minimum length sequence needed to accommodate option $o$ occurring $t$ times, i.e., $\mu_o(t) = p_o((t-1)$ \textbf{div} $q_o) + ((t-1)$ \textbf{mod} $q_o) + 1$ \citep{perron2004combining}.  


\added[]{
The features 1--12 have already been used in different ways to create current benchmark problem instances \citep{gent1999csplib,gravel04,perron04}. Hence, we use them as part of our feature set. Moreover, we have seen in \citet{thiruvady2020large} that the LCM of $q$ values typically leads to effective subsequences that can be solved efficiently, and hence we make use of this feature 13. Finally, all known problem instances have been obtained by using uniform distributions to assign cars to classes, whereas allowing variances in these assignments could impact the hardness of a problem. Hence, we also include feature 14 in our feature set.      
}

These features will be used to perform instance space analysis and also used as inputs into machine learning models to build algorithm selection models. Hence, they form the basis of this study. Moreover, we will systematically generate new problem instances by varying these features.

\subsection{Instances}
\label{insts}

We use a total of 247 problem instances in this study, including 121 existing benchmark instances and 126 newly generated instances. 


The existing benchmark instances we have chosen consist of the nine instances available from CSPLIB \citep{gent1999csplib}, 82 from \citep{perron2004combining}, and 30 from \citep{gravel2005review}. These benchmark instances were previously considered ``hard'' to solve, though, the last decade and a half has seen rapid development of commercial solvers and extensive research in heuristic/meta-heuristic methods, a number of these problem instances can now be solved efficiently. Hence, we also generate new problem instances with various problem characteristics for stress testing of algorithm performance. 


The new problem instances are generated as follows. To generate an instance specified with $n$ cars, $o$ options, and $k$ car classes, we generate $o$-dimensional option vectors, by including option $i$ with probability $p_i/q_i$, until we reach $k$ unique vectors. We then assign each car to a car class with uniform probability. Each of the following fourteen sets contain nine instances, three for each size of 100, 300, and 500 cars, for an overall total of 126 new instances:
\begin{enumerate}
  \item \textsf{noBhiU}: no bias and high utilisation. This is the ``default'' set with $o = 5$, $k = 25$, $1 \leqslant p \leqslant 3$, and $1 \leqslant q - p \leqslant 2$. The car classes are selected uniformly on a per-car basis, and utilisation rate for all options is between $90-100\%$.
  \item \textsf{negBhiU}: negative bias and high utilisation. It is the same as \textsf{noBhiU}, except that car classes with fewer options receive more cars. In particular, when a car is assigned to a given class $k$ containing $o_k$ options, we randomly reassign it with probability $o_k/o$. 
  \item \textsf{posBhiU}: positive bias and high utilisation. It is the same as \textsf{noBhiU}, except that car classes with more options receive more cars. In particular, when a car is a assigned to a given class $k$ containing $o_k$ options, we randomly reassign it with probability $1 - o_k/o$. 
  \item \textsf{hiPQhiU}: high $p/q$ ratios and high utilisation. It is the same as \textsf{noBhiU}, except that $2 \leqslant p \leqslant 4$.
  \item \textsf{hiPQmedU}: high $p/q$ ratios and medium utilisation. It is the same as \textsf{hiPQhiU}, except that all option utilisation rates are between $70-80\%$.
  \item \textsf{negBhiPQloU}: negative bias, high $p/q$ ratios, and low utilisation. It is the same as \textsf{hiPQhiU}, except that all option utilisation rates are between $50-60\%$, and negative bias is added. It is difficult to generate such low-utilisation instances using $p/q$ option probabilities without negative bias.
  \item \textsf{loPQ}: low $p/q$ ratios. The settings for \textsf{loPQ} are: $o = 8$, $k = 20$, $1 \leqslant p \leqslant 2$, $ 2 \leqslant q - p \leqslant 3$, and no lower bound on utilisation. It is difficult to generate low-$p/q$-ratio instances with $o = 5$ and $k = 25$.
  \item \textsf{negBfixedPQ}: negative bias and fixed $p$ and $q$. The parameter $p$ is fixed to $(3, 2, 1, 2, 1)$; $q$ is fixed to $(4, 3, 4, 5, 2)$; and there is no lower bound on utilisation. It is difficult to generate these instances (presumably again because of the low $p/q$ ratio included, namely $1/4$) without negative bias.
  \item \textsf{randN}: random configuration of car classes. It is the same as \textsf{noBhiU}, except that we select uniformly from the $\binom{n+k-1}{k-1}$ possible assignments of cars to car classes, instead of uniform distribution of cars to classes. Doing so increases the variance in car class populations.
  \item \textsf{RloU}: random car classes and low utilisation (with option utilisation rate only restricted to $>50\%$).
  \item \textsf{RanyU}: random car classes, no utilisation restrictions, and 10 options.
  \item \textsf{RnegBvloU}: random car classes, negative bias, and very low utilisation (with option utilisation rate $<60\%$).
  \item \textsf{RnegBhiPQvloU}: random car classes, negative bias, high $p/q$ ratios, and very low utilisation (with option utilisation rate $<60\%$).
  \item \textsf{RnegBloPQvloU}: random car classes, negative bias, low $p/q$ ratios, and very low utilisation (with option utilisation rate $<60\%$).
\end{enumerate}
While the above sets of parameter choices may seem  arbitrary, they were selected after significant experimentation to provide a good coverage of the instance space, as will be shown in Section~\ref{sec:visualInst}.

\section{Experimental Results}
\label{sec::experiments}
This section presents the experimental results of the instance space analysis and algorithm selection for the car sequencing problem. The solution algorithms were implemented in C++ and compiled with GCC-5.2.0, \added[]{and the algorithm selection models were implemented in Python with the scikit-learn library \citep{scikit-learn}.} The experiments were conducted on the Monash University campus cluster, MonARCH. The source code and problem instances will be made publicly available when the paper is published.

\subsection{Visualizing Instances}\label{sec:visualInst}
We create a two-dimensional instance space to visualize how well the newly generated problem instances complement the existing instances available in the literature. Each instance is represented by a vector of fourteen features listed in Section~\ref{features}. For features that are not in the range of $0$ and $1$ (i.e., number of cars, number of options, number of car classes, and lowest common multiple of $q$ values), we normalize them by their maximum value. For better visualization, we select eight features that are mostly correlated with the algorithm performances using information theoretic feature selection methods \citep{sun2020revisiting}. The eight features selected are the option utilisation (minimum, average, standard deviation), number of car classes, number of options, average number of options per car class and $p/q$ ratio (average and maximum). 

\begin{figure}[!t]
    \centering
    \resizebox{\textwidth}{!}{ 
    \includegraphics[scale=0.30]{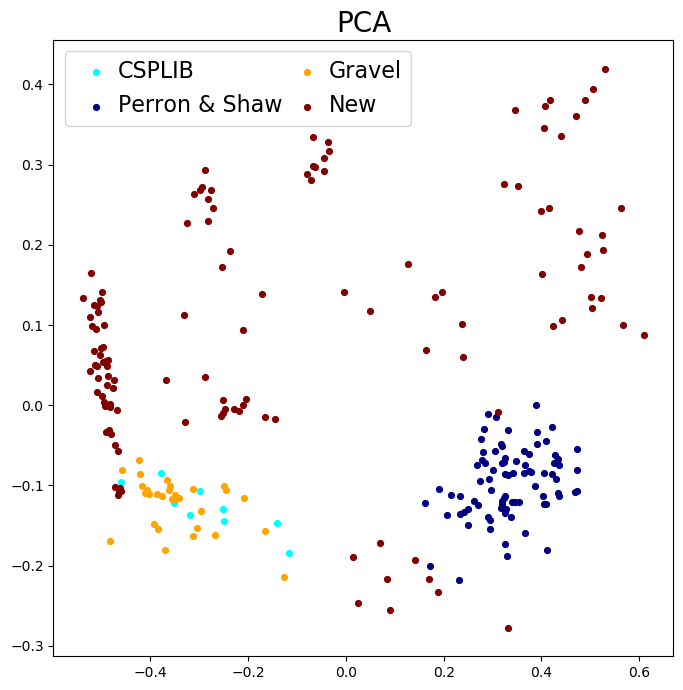}
    \includegraphics[scale=0.3]{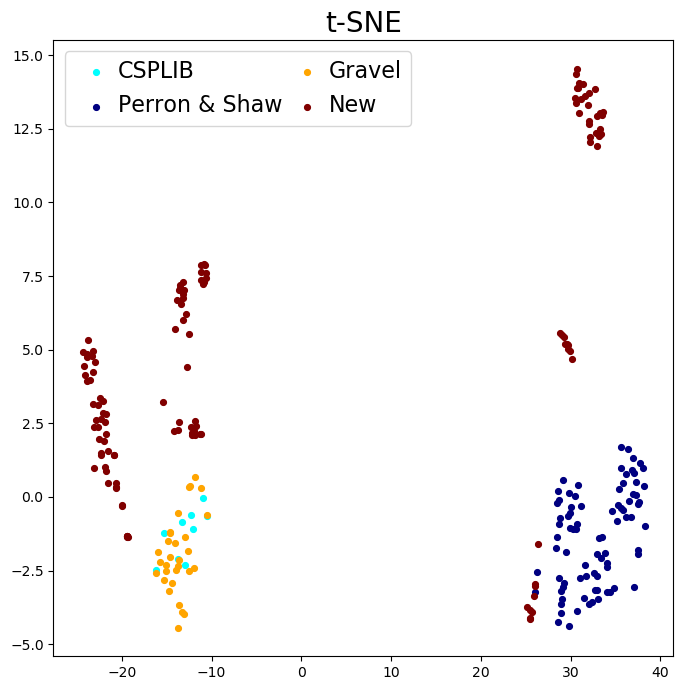}
    \includegraphics[scale=0.3]{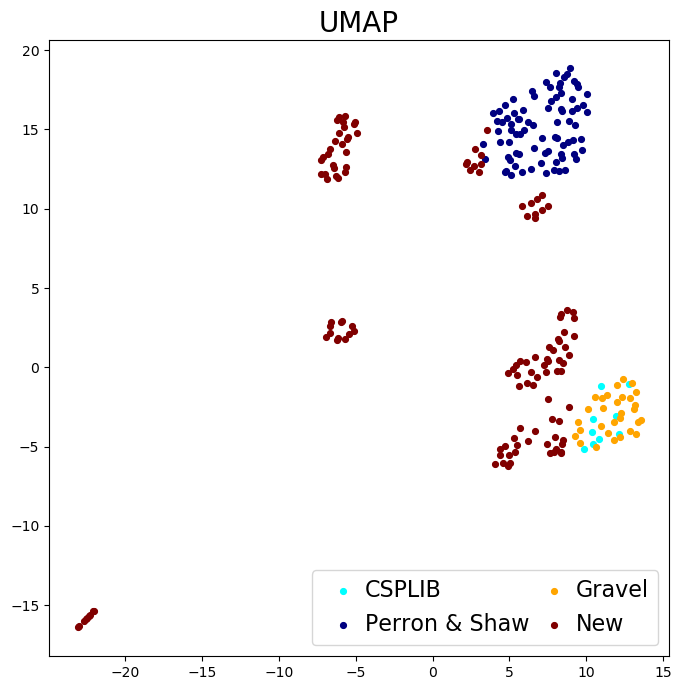}
    }
    \caption{\added[]{The two-dimensional instance spaces created by PCA, t-SNE and UMAP. Here each subfigure, a dot represents a problem instance, and its color indicates which category the instance belongs to.}}
    \label{fig:instance distribution}
\end{figure}

\begin{figure}[!t]
\centering
\resizebox{0.40\textwidth}{!}{
\begin{tikzpicture}
	\begin{axis} [box plot width=0.20em, xtick=data,  ylabel = explained variance percentage (\%),  xlabel = number of components selected, height=0.51\textwidth,width=0.70\textwidth, grid style={line width=.1pt, draw=gray!10},major grid style={line width=.2pt,draw=gray!30}, xmajorgrids=true, ymajorgrids=true,  major tick length=0.05cm, minor tick length=0.0cm, legend style={at={(0.25,0.80)},anchor=west,font=\scriptsize,draw=none}]
	\addplot[color=marine, mark=*, mark size = 1.0, line width=0.40mm] coordinates {
	(1, 0.7456066715810046*100)
    (2, 0.8763205054861547*100)
    (3, 0.9376304867234331*100)
    (4, 0.9734269746019238*100)
    (5, 0.9864111292141873*100)
    (6, 0.9969198636884883*100)
    (7, 0.9994071912662297*100)
    (8, 1.0*100)
	}; 
	\end{axis}
\end{tikzpicture}
}
\caption{\added[]{The percentage of variance explained versus the number of components selected in PCA. The first two components capture 88\% of the data variance.}}
\label{fig:PCA variance}
\end{figure}
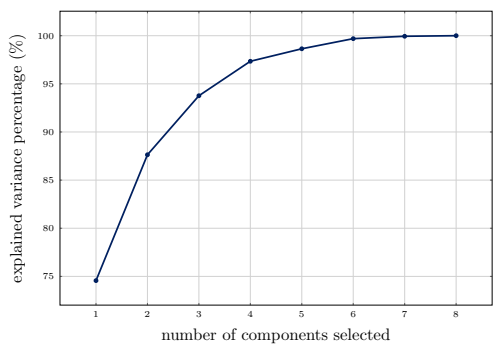

We then use dimensionality reduction techniques to project the eight-dimensional feature vectors onto a two-dimensional plane. \added[]{The techniques used are principal component analysis (PCA), t-distributed stochastic neighbor embedding (t-SNE) \citep{van2008visualizing} and uniform manifold approximation and projection (UMAP) \citep{mcinnes2018umap}. The created two-dimensional instance spaces are shown in Fig.~\ref{fig:instance distribution}. The instances are more evenly spread in the space created by PCA compared to the other two, and therefore we will only use PCA for visualization in the rest of the paper. The percentage of variance explained versus the number of components selected in PCA is shown in Fig.~\ref{fig:PCA variance}. The first two principal components explain 88\% of the data variance,} and the corresponding transformation are listed as follows:
\begin{gather}
 Z =
  \begin{bmatrix}
    0.46417697  & 0.13118315 \\
    -0.13786325 & 0.40042554 \\
    -0.71876354 & 0.03250308 \\
    -0.36351993 & -0.44927092 \\
    0.22506848 & -0.33242684 \\
    -0.25113443 & 0.21863308 \\
    -0.04415605 & 0.59942745 \\
    0.03303967 &  0.31926206 
   \end{bmatrix}^T
   \begin{bmatrix}
    \text{num-options} \\
    \text{num-classes} \\
    \text{min-utilisation} \\
    \text{ave-utilisation} \\
    \text{std-utilisation} \\
    \text{ave-options} \\
    \text{ave-pq-ratio} \\
    \text{max-pq-ratio} 
   \end{bmatrix}.
\end{gather}

\begin{figure}[!t]
	\centering
	\resizebox{\textwidth}{!}{ 
		\includegraphics[scale=0.30]{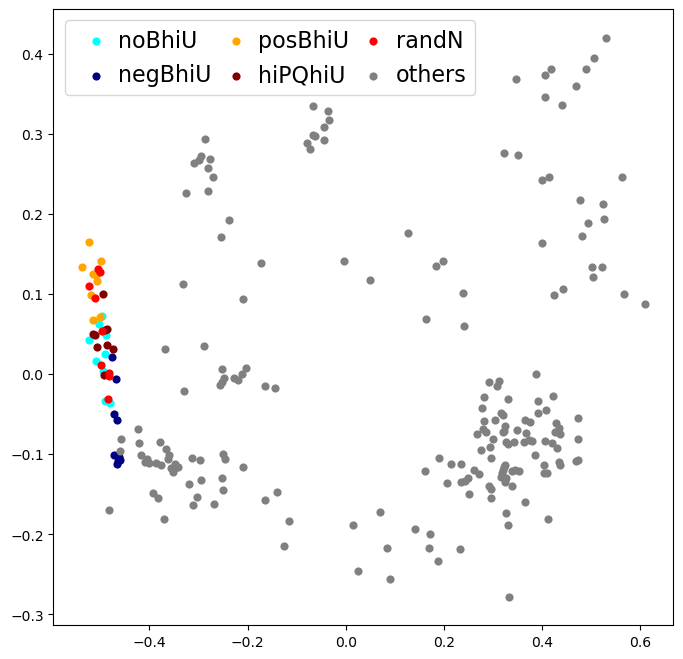}
		\includegraphics[scale=0.30]{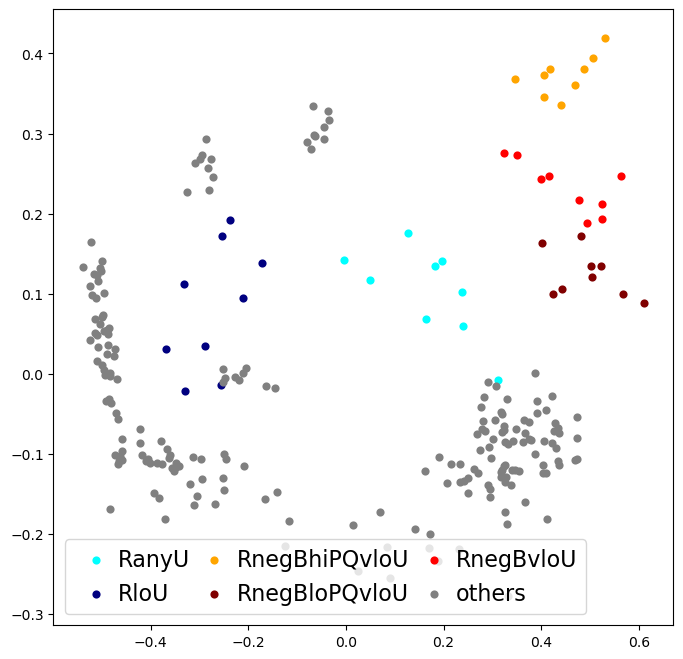}
		\includegraphics[scale=0.3]{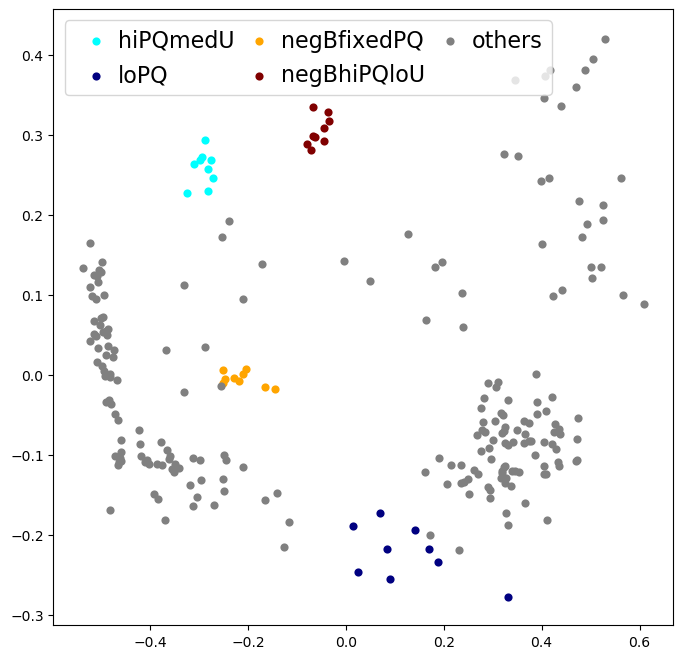}
	}
	\caption{The distribution of the newly generated instances. The axes of the figures are the first two principal components from PCA, and each dot represents a problem instance.}
	\label{fig: new instance distribution}
\end{figure}

We can observe that the existing problem instances, i.e., CSPLIB, Gravel, and Perron \& Shaw, mainly lie at the bottom of the two-dimensional instance space created via PCA. The newly generated instances complement the existing instances, and they together fill the instance space well. We also visualize where each category of the newly generated instances lies in the two-dimensional instance space in Fig. \ref{fig: new instance distribution}. 

\begin{figure}[!t]
    \centering
    \includegraphics[scale=0.36]{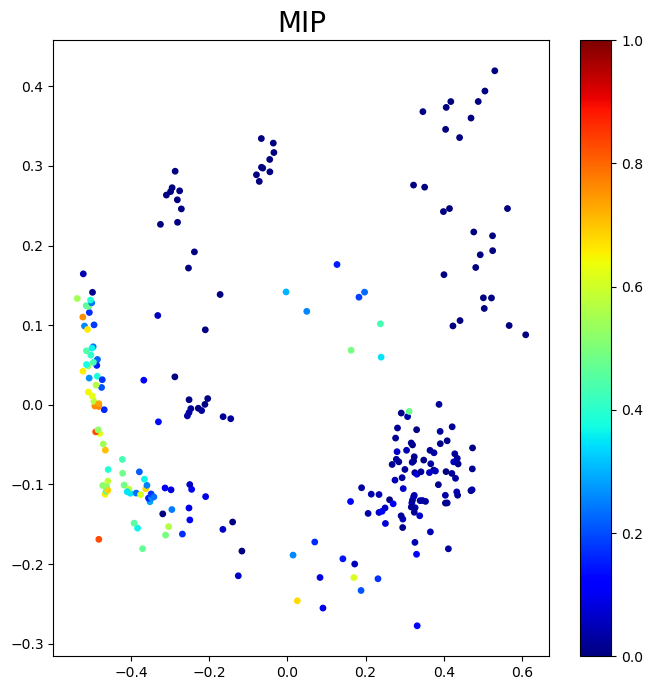}
    \includegraphics[scale=0.36]{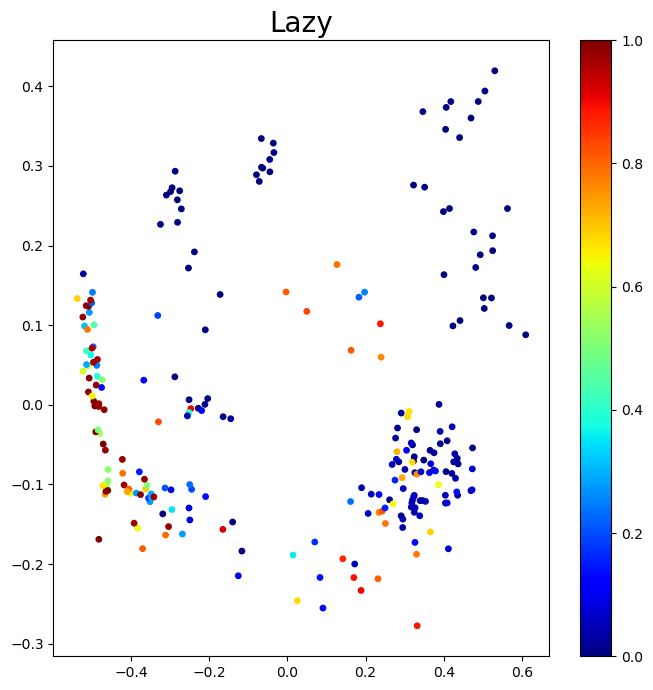}

    \includegraphics[scale=0.36]{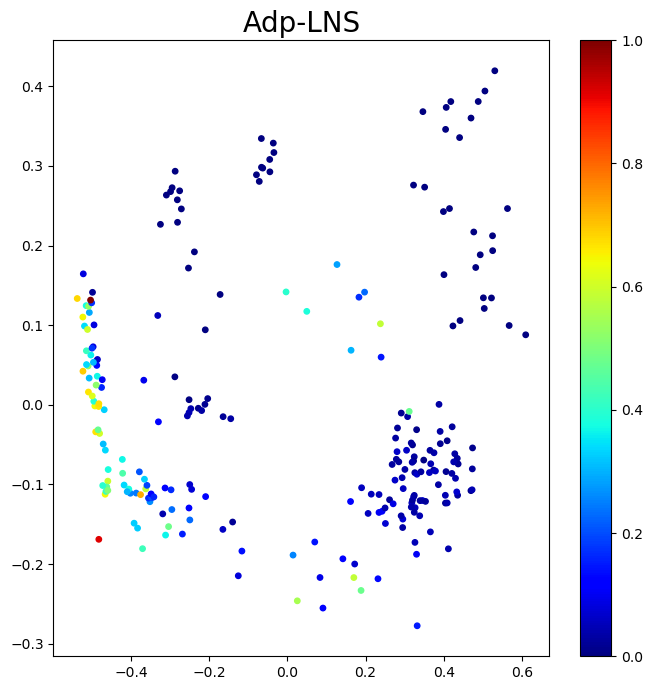}
    \includegraphics[scale=0.36]{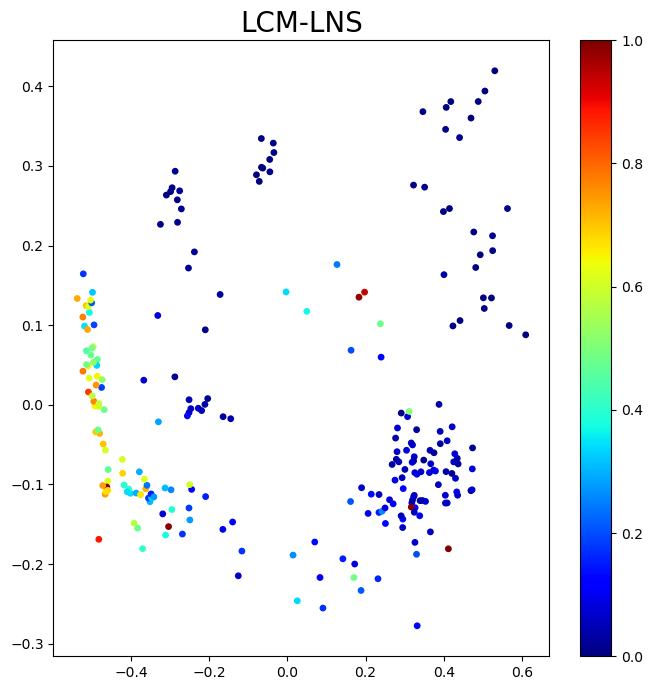}
  
    \includegraphics[scale=0.36]{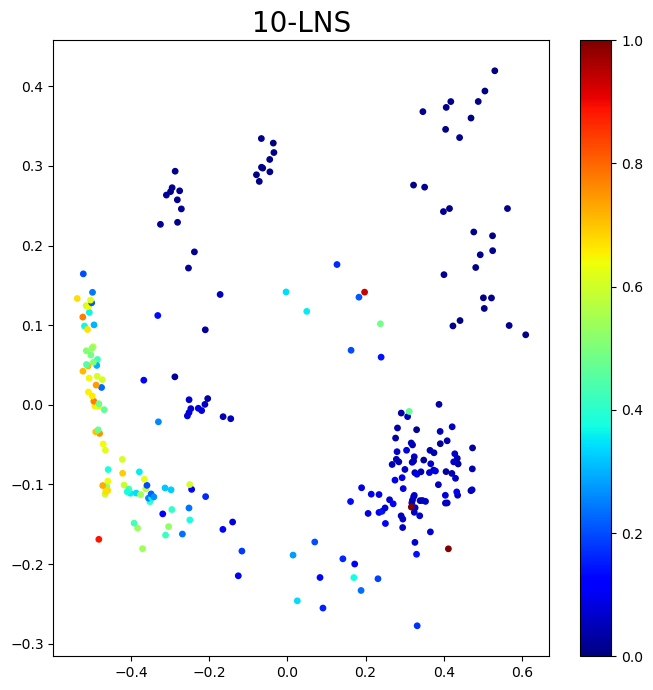}
    \includegraphics[scale=0.36]{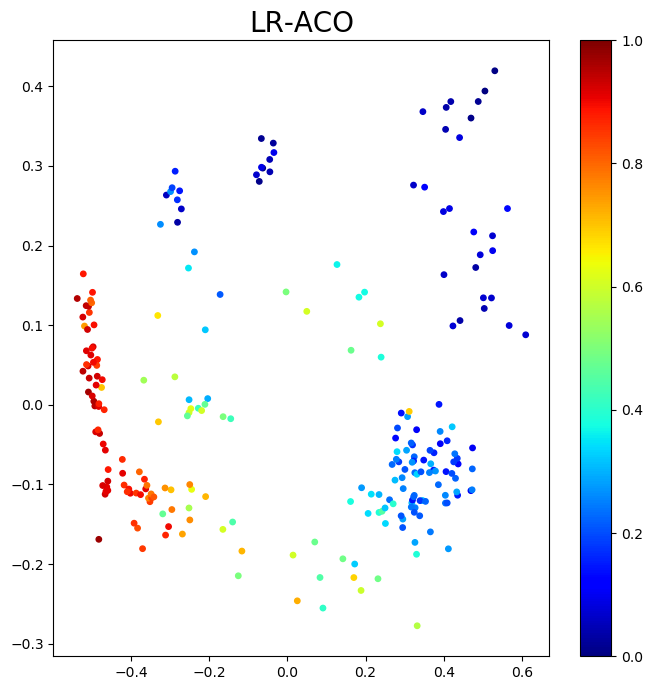}
    \caption{The optimality gap generated by each algorithm for the problem instances tested. The axes of the figures are the first two  principal components of the feature vectors; each dot represents a problem instance in the space spanned by the first two  principal components; and the color of the dots indicates the optimality gap (red is larger).}
    \label{fig:optimality gap}
\end{figure}

\subsection{Visualizing Algorithm Performance}
We test the performance of six algorithms on the 247 problem instances. The algorithms are LR-ACO, 10-LNS (LNS with a fixed window size of 10),  LCM-LNS (LNS with a window size computed by the lowest common multiple of $q$ values), and MIP (solved by Gurobi), as well as two new algorithms i.e., Adp-LNS (the adaptive version of LNS) and  MIP with lazy constraints (solved by Gurobi). The parameter setting for LR-ACO is consistent with the original paper \citep{Thiruvady2014a}. Gurobi 6.5.1 is used to solve the MIPs with 4 cores and 10GB of memory given. The parameter settings for Gurobi are tuned by hand: barrier iteration limit = 30; MIPgap = 0.005; Presolve = 1 (conservative); barrier ordering = nested dissection; cuts = -1 (aggressive); cut passes = 200 (large); barrier convergence tolerance = 0.0001; presparsify = 1 (turn on always); pre pass limit = 2. The time limit for each algorithm to solve an instance is set to one hour. The objective value ($y_b$) generated by an algorithm for an instance is compared against the best known lower bound ($y_{lb}$) for that instance to compute an optimality gap $(y_b-y_{lb})/y_b \in [0,1)$. The gap generated by an algorithm for an instance indicates how well the algorithm performs on that instance. 

The performance of each algorithm is visualized in the two-dimensional instance space in Fig. \ref{fig:optimality gap}, where each dot represents an instance and its color indicates the optimality gap generated by the corresponding algorithm for that instance. We observe that most of the hard instances lie in the bottom left region of the instance space, which corresponds to five categories of newly generated instances (as shown in the first graph of Fig. \ref{fig: new instance distribution}), as well as part of the CSPLIB and Gravel instances. Some of the newly generated \textsf{RanyU} and \textsf{loPQ} instances are also hard. The Perron \& Shaw instances were considered to be among the hardest problem instances available in the literature \citep{perron2004combining}, but we observe that these instances are relatively easy for the six algorithms to solve. The LR-ACO algorithm is clearly outperformed by other algorithms, thus we will leave LR-ACO out in the subsequent analysis. 

\begin{figure}[!t]
    \centering
    \includegraphics[scale=0.28]{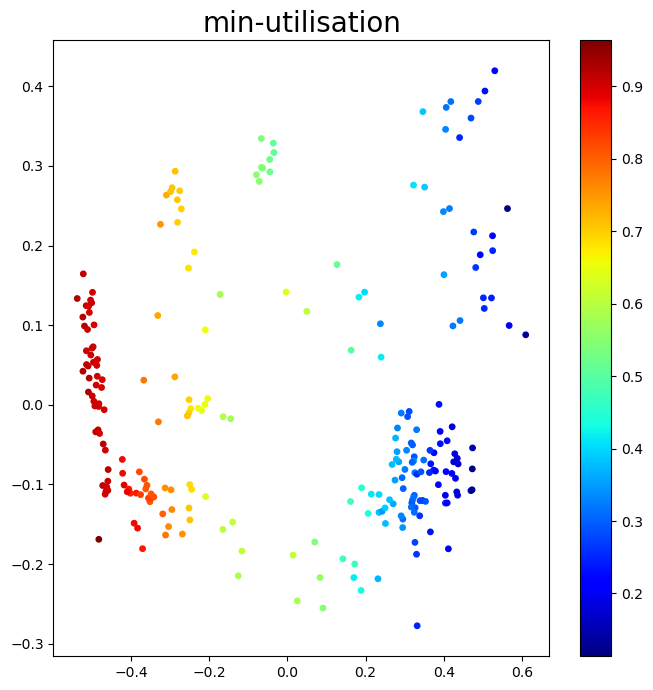}
    \includegraphics[scale=0.28]{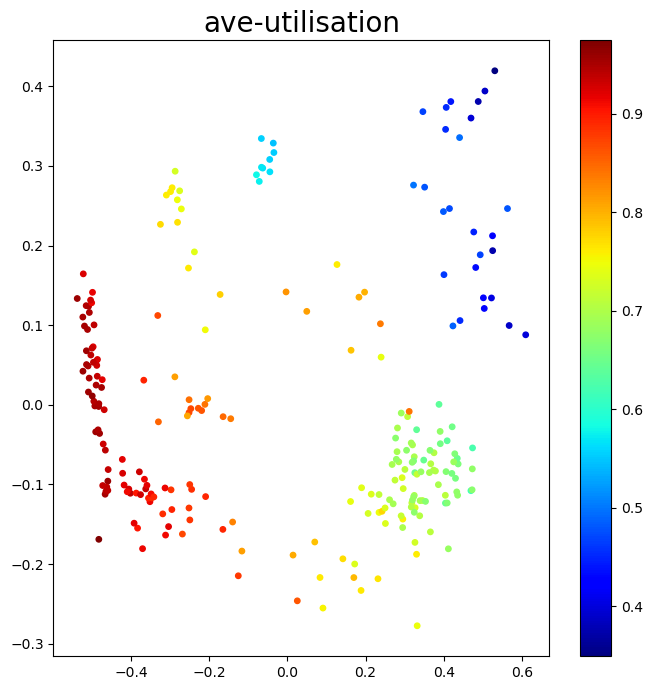}
   
    \includegraphics[scale=0.28]{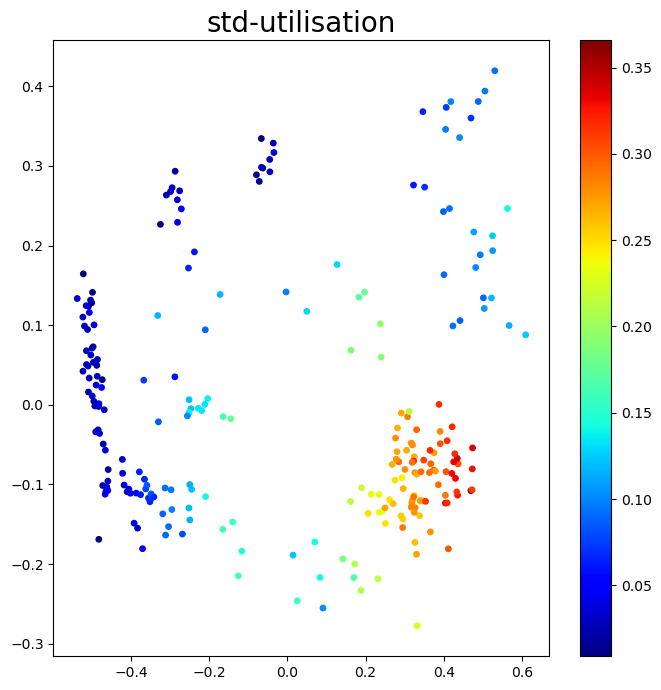}
    \includegraphics[scale=0.28]{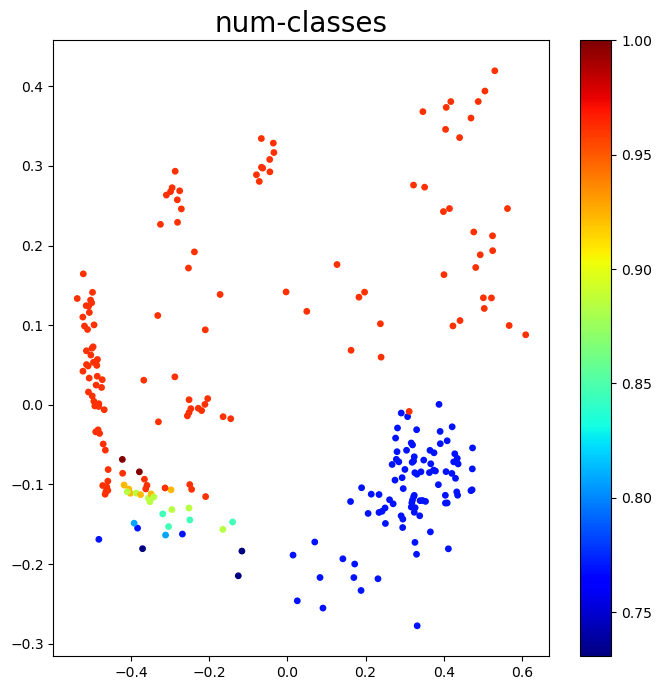}
    
    \includegraphics[scale=0.28]{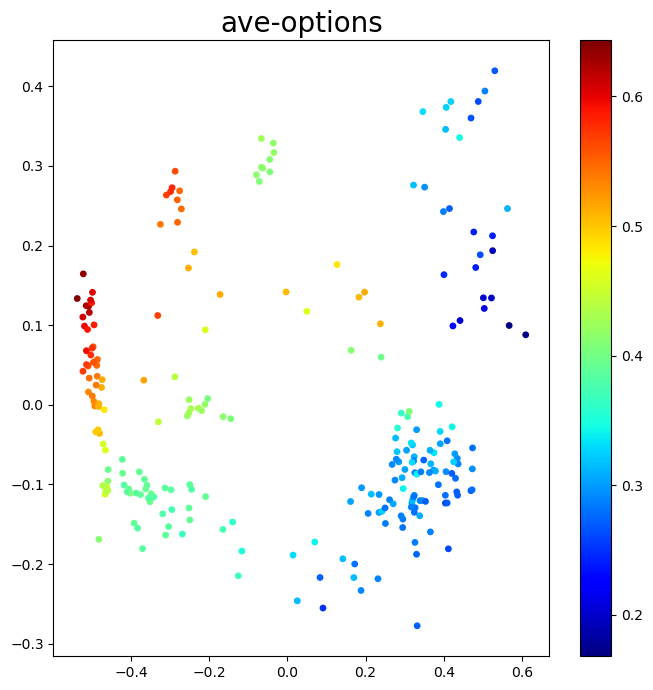}
    \includegraphics[scale=0.28]{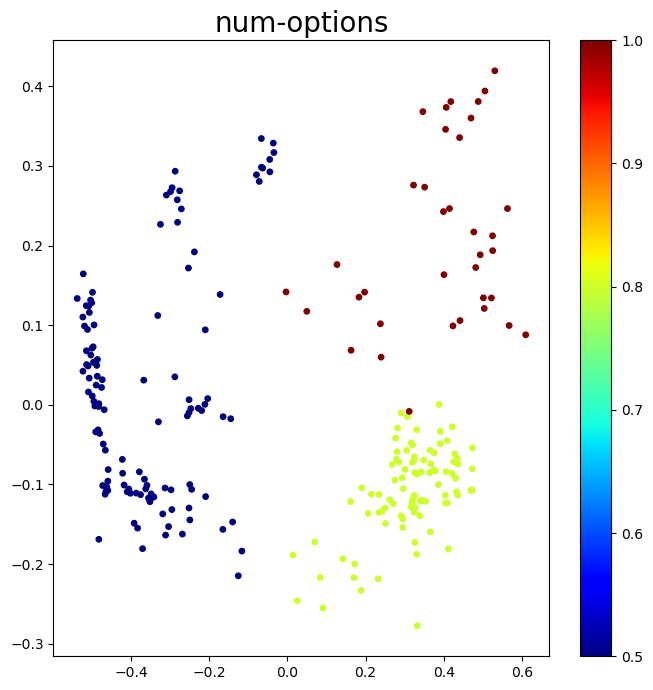}
    
    \includegraphics[scale=0.28]{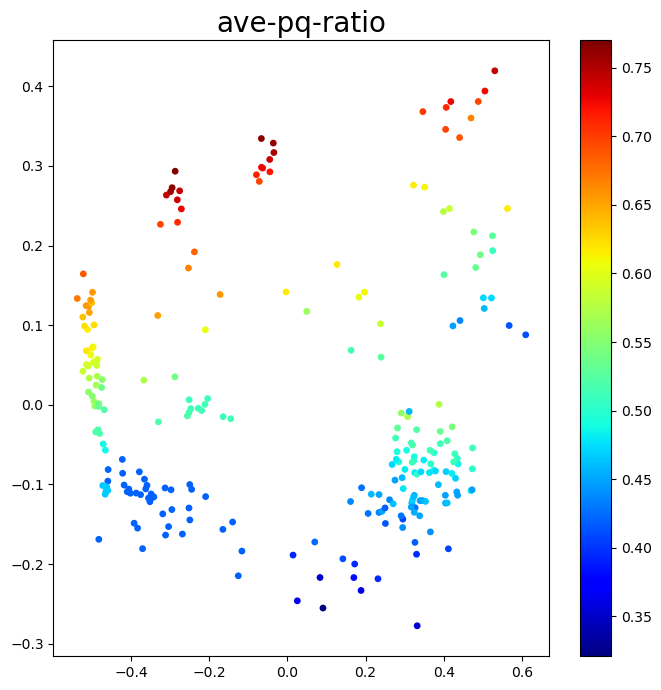}
    \includegraphics[scale=0.28]{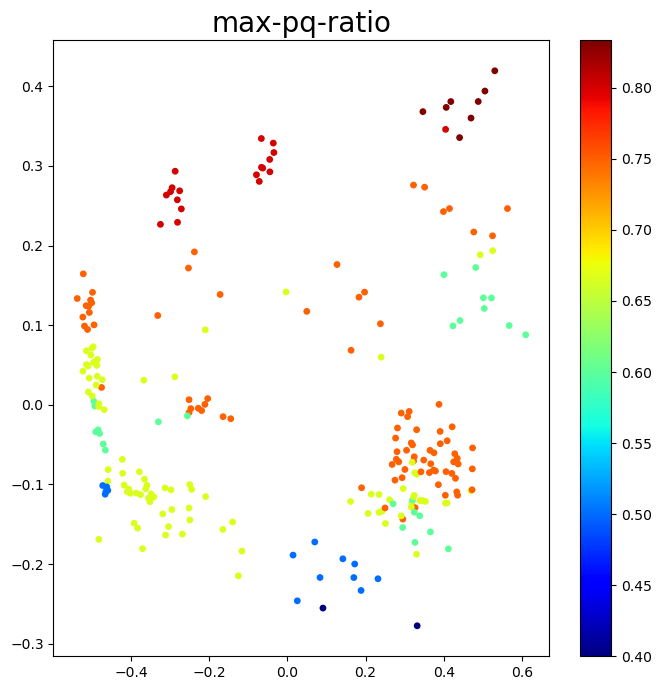}
    \caption{The distribution of each feature in the two-dimensional instance space created via PCA. Each dot represents a problem instance, and its color indicates the value of the corresponding feature of that instance.}
    \label{fig:feature distribution}
\end{figure}

\subsection{Visualizing Features}
We visualize the eight features in the two-dimensional instance space in Fig.~\ref{fig:feature distribution}, where each dot represents an instance and its color represents the value of corresponding feature. The aim is to identify which features make an instance hard to solve. From Fig.~\ref{fig:optimality gap} and Fig.~\ref{fig:feature distribution}, we can observe that the minimum utilisation (or average utilisation) is highly positively correlated with the optimality gap generated by the algorithms, meaning that high utilisation makes an instance hard to solve. This result is not surprising since higher utilisation effectively means that satisfying the subsequence constraints is a lot more complex, and as a consequence, the modulation of the use of option is also more complex. Instances with a large average number of options per car class are mostly hard to solve. For the same reasons previously discussed, this result is also expected. However, the number of options is negatively correlated with the optimality gap, indicating that a large number of options makes an instance easy to solve. This is an interesting outcome which shows that a large number of options allows less room for car classes to be sequenced in different positions. There is no significant correlation between the algorithm performance and other three features, i.e., number of car classes, average p/q ratio and maximum p/q ratio. 

\subsection{Algorithm Selection}
\added[]{We build our algorithm selection model based on machine learning using the full set of fourteen features (Section~\ref{features}) to select the best algorithm for solving a given instance.} We perform four comparisons between the algorithms: 1) MIP vs Lazy, 2) Adp-LNS vs LCM-LNS, 3) Adp-LNS vs 10-LNS, and 4) MIP vs Adp-LNS. For each pair of algorithms (denoted as A and B), we divide the instances into three classes: 1) algorithm A performs better; 2) algorithm B performs better; or 3) algorithms A and B perform equally well if the difference between the optimality gaps generated are less than 0.005. It then becomes the standard three-class classification problem, and any classification algorithm can be used for this task. We have tested support vector machine (SVM) with a radial basis function kernel, k nearest neighbour (KNN), decision tree (DT) \added[]{and logistic regression (LR) models. We have tuned the regularization parameter ($C$) of SVM and LR, number of neighbors ($k$) of KNN, and maximum tree depth ($d$) of DT, with the ten-fold cross-validation accuracy reported in Table~\ref{Table: parameter tunning}. The best parameter values identified are $C=1000$ for SVM, $k=3$ for KNN, $d=2$ for DT and $C=10$ for LR. The other parameter settings are consistent with the default of the scikit-learn library.}


\begin{table}[!t]
\centering
\caption{The ten-fold cross-validation accuracy of algorithm selection generated by the classifiers with different parameter settings.}
\label{Table: parameter tunning}
\resizebox{\textwidth}{!}{ 
\begin{tabular}{lllll|llll|llll|llll}
\toprule
\multirow{2}{*}{Algorithm Pairs} & \multicolumn{4}{c|}{SVM ($C$) }     &  \multicolumn{4}{c|}{KNN ($k$)}  &    \multicolumn{4}{c|}{DT ($d$)}  & \multicolumn{4}{c}{LR ($C$)} \\\cline{2-5}\cline{6-9}\cline{10-13}\cline{14-17}
& 1 & 10 & 100 & 1000 & 2 & 3 & 4 & 5 & 2 & 3 & 4 & 5 & 1 & 10 & 100 & 1000 \\\midrule
MIP vs Adp-LNS  & 0.43 & 0.64  & 0.69   & 0.69  & 0.52 & 0.58 & 0.55 & 0.53 & 0.72 & 0.69 & 0.70 & 0.67 & 0.69 & 0.71 & 0.69  & 0.68   \\
MIP vs Lazy     & 0.57 & 0.75  & 0.79   & 0.83  & 0.75 & 0.75 & 0.72 & 0.71 & 0.81 & 0.84 & 0.81 & 0.82 & 0.78 & 0.82 & 0.83  & 0.83   \\
LNS (Adp vs LCM) & 0.76 & 0.77  & 0.81   & 0.83 & 0.74 & 0.74 & 0.75 & 0.76 & 0.80 & 0.79 & 0.82 & 0.82 & 0.78 & 0.82 & 0.82  & 0.82   \\
LNS (Adp vs 10)  & 0.89 & 0.89  & 0.90   & 0.90 & 0.89 & 0.88 & 0.89 & 0.89  & 0.87 & 0.87 & 0.87 & 0.87  & 0.89 & 0.91 & 0.90  & 0.90   \\\bottomrule  
\end{tabular}
}
\end{table}

\begin{table}[!t]
\centering
\caption{The average ten-fold cross-validation accuracy of algorithm selection.}
\label{table: selection accuracy}
\setlength{\tabcolsep}{8pt}
\resizebox{\textwidth}{!}{ 
\begin{tabular}{llllll}
\toprule
Algorithm Pairs  & SVM ($C=10^3$)  & KNN ($k=3$)  & DT ($d=2$)  & LR ($C=10$)  & ZeroR \\\midrule
MIP vs Adp-LNS   & 0.694 & 0.582 & \bf{0.716} & 0.706 & 0.381 \\
MIP vs Lazy      & \bf{0.833} & 0.741 & 0.802 & 0.818 & 0.567 \\
LNS (Adp vs LCM) & \bf{0.833} & 0.739 & 0.799 & 0.822 & 0.757 \\
LNS (Adp vs 10)  & 0.895 & 0.872 & 0.872 & \bf{0.906} & 0.891 \\
\bottomrule
\end{tabular}
}
\end{table}


\begin{figure}[!t]
    \centering
    \includegraphics[scale=0.35]{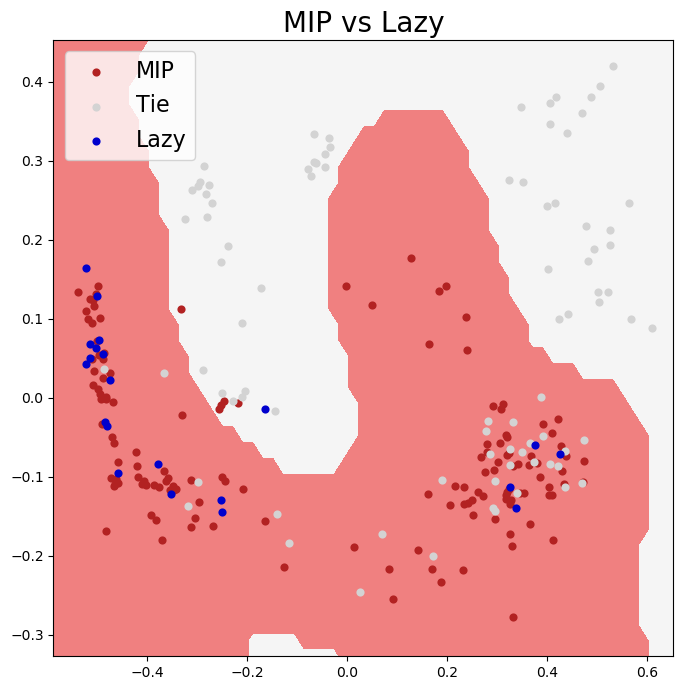}
    \includegraphics[scale=0.35]{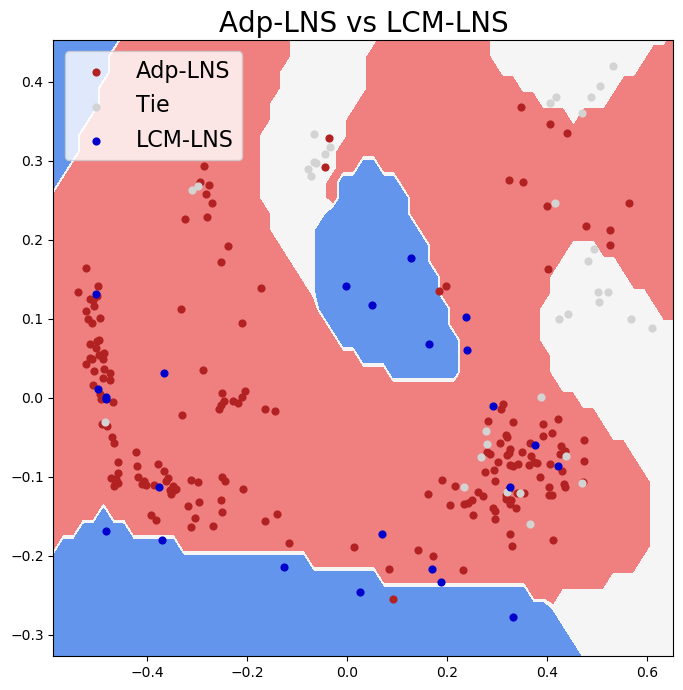}
    
    \includegraphics[scale=0.35]{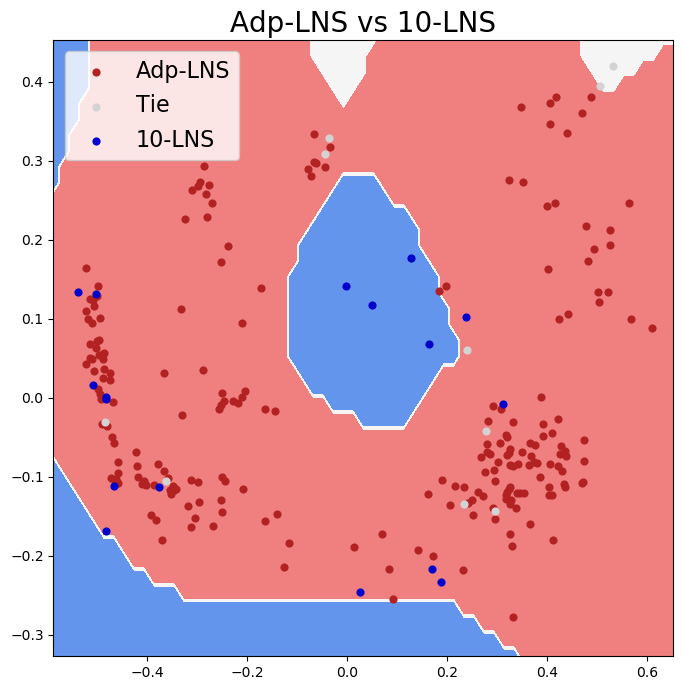}
    \includegraphics[scale=0.35]{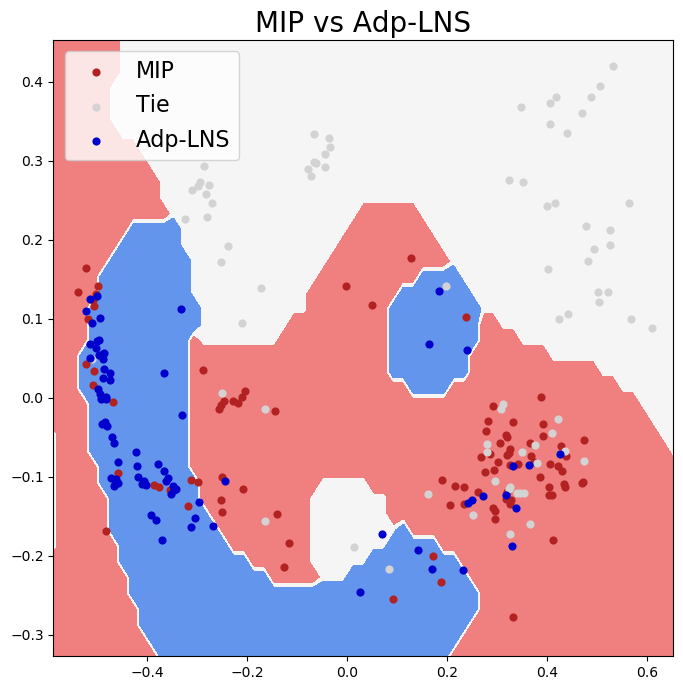}
    \caption{Pairwise comparison between the algorithms in terms of the optimality gap generated. The axes of the figures are the first two  principal components of the feature vectors. Each dot represents a problem instance, and the color of the dots indicates which algorithm performs better. The color of a region indicates which algorithm performs better in that region predicted by SVM.}
    \label{fig:pairwise comparision}
\end{figure}

The average ten-fold cross-validation accuracy of algorithm selection over 100 runs using each machine learning model is reported in Table \ref{table: selection accuracy}. The results show that the machine learning models achieve a reasonably good accuracy for the algorithm selection tasks, comparing to the ZeroR method that simply predicts the majority class. KNN is clearly outperformed by the other three machine learning models. 

\begin{figure}[!t]
    \centering
    \resizebox{0.6\textwidth}{!}{
    \begin{tikzpicture}[shorten >=1pt,auto, node distance=1cm, thick,main node/.style={color =marine, circle,draw,font=\sffamily,inner sep = 5pt}]
    \node[main node, rectangle] (1) at (0,0) {\small max-utilisation $\le$ 0.885};
    \node[main node, rectangle] (2) at (-3,-2) {\small ave-pq-ratio $\le$ 0.366 };
    \node[main node, rectangle] (3) at (3,-2) {\small min-utilisation $\le$ 0.608 };
    \node[main node, rectangle, minimum width=1.5cm] (4) at (-4.5,-4) {\small MIP };
    \node[main node, rectangle, minimum width=1.5cm] (5) at (-1.5,-4) {\small Tie };
    \node[main node, rectangle, minimum width=1.5cm] (6) at (1.5,-4) {\small MIP };
    \node[main node, rectangle, minimum width=1.5cm] (7) at (4.5,-4) {\small Adp-LNS };
    \path[every node/.style={color = marine,font=\sffamily\small}]
    (1) [color = marine, ->, sloped] edge node {\small Y} (2)  (1) edge node {\small N} (3) (2) edge node {\small Y} (4) (2) edge node {\small N} (5) (3) edge node {\small Y} (6) (3) edge node {\small N} (7);
    \end{tikzpicture}
    }
    \caption{\added[]{The decision tree for selecting the best algorithm between MIP and Adp-LNS to solve a problem instance.}}
    \label{fig:DT tree}
\end{figure}
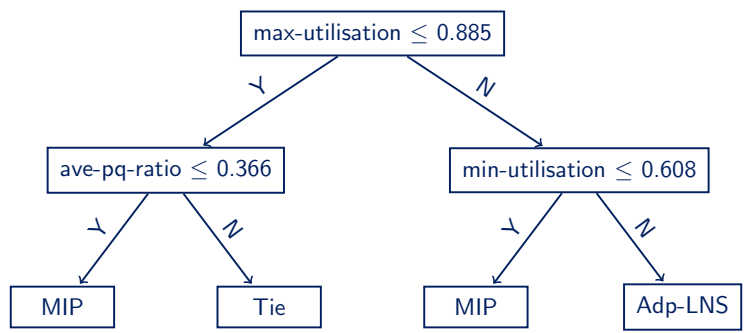

To gain a deeper insight into the algorithm performance, we visualize the decision boundary identified by SVM in the two-dimensional instance space in Fig. \ref{fig:pairwise comparision}. The color of a dot represents which algorithm performs better on that instance and the color of a region indicates which algorithm performs better for that region based on the SVM prediction. We can clearly see that MIP performs equally well as Lazy on easy instances, and outperforms Lazy on other instances. The Adp-LNS algorithm outperforms the other two LNS algorithms on most of the instances except \textsf{RanyU} instances. We then take the best performing exact algorithm, MIP, and compare it against the best performing local search method, Adp-LNS. Interestingly, by comparing Fig. \ref{fig:pairwise comparision} with Fig. \ref{fig:optimality gap}, we can see a clear pattern that (1) Adp-LNS outperforms MIP on hard instances; (2) MIP outperforms Adp-LNS on medium-hard instances; and (3) they do equally well on easy instances. \added[]{Fig.~\ref{fig:DT tree} presents the decision tree learned by the DT classifier to select the best algorithm between MIP and Adp-LNS for solving a problem instance. These results provide great insights into the algorithms' performance, and a simple guideline on which algorithm should be selected to solve a given car sequencing problem instance in practice.}  

\added[]{Finally, we report the AUC (area under the receiver operating characteristic curve), precision, recall and F1-score generated by each machine learning model for the algorithm selection tasks in Table~\ref{tab::auc}. The machine learning models achieve a reasonably good performance in classifying majority class, while they may have trouble in predicting minority class when the data for the classification task is extremely unbalanced. For example, the machine learning models do not perform well in classifying the instances where 10-LNS outperforms Adp-LNS. However, this result is acceptable in practice because Adp-LNS generally performs much better than 10-LNS, and simply selecting Adp-LNS to solve all the instances would generate the best results for 94\% of the instances when compared to 10-LNS.} 

\begin{table}[!t]
\centering
\caption{\added[]{The AUC, precision, recall and F1-score generated by each machine learning model for the algorithm selection tasks. For a pair of algorithms A and B, the classification results for each of the three classes are reported: (win) A better than B, (tie) A equal to B, and (loss) A worse than B.}}
\label{tab::auc}
\resizebox{\textwidth}{!}{ 
\begin{tabular}{ll|tsw|tsw|tsw|tsw}
\toprule
\multirow{2}{*}{Algorithm Pairs}  & \multirow{2}{*}{Metrics} &      \multicolumn{3}{c|}{SVM ($C=10^3$) }     &  \multicolumn{3}{c|}{KNN ($k=3$)}  &    \multicolumn{3}{c|}{DT ($d=2$)}  & \multicolumn{3}{c}{LR ($C=10$)} \\\cline{3-5}\cline{6-8}\cline{9-11}\cline{12-14}
& & win & tie & loss & win & tie & loss & win & tie & loss & win & tie & loss \\\midrule
\multirow{4}{*}{MIP vs Adp-LNS}  & AUC  & 0.73    & 0.85    & 0.83    & 0.66    & 0.67    & 0.70    & 0.69    & 0.85    & 0.80    & \bf{0.77}    & \bf{0.86}    & \bf{0.86}    \\
& precision & \bf{0.59} & 0.89 & 0.65 & 0.49 & 0.59 & 0.53 & \bf{0.59} & \bf{0.96} & \bf{0.73} & \bf{0.59} & 0.95 & 0.71 \\
& recall    & 0.73 & 0.61 & 0.66 & 0.69 & 0.33 & 0.48 & \bf{0.77} & \bf{0.62} & 0.72 & 0.74 & \bf{0.62} & \bf{0.73} \\
& F1-score    & 0.65 & 0.72 & 0.65 & 0.57 & 0.41 & 0.50 & \bf{0.67} & \bf{0.75} & \bf{0.72} & 0.65 & \bf{0.75} & 0.71 \\\midrule
\multirow{4}{*}{MIP vs Lazy}    & AUC          & 0.92    & 0.92    & 0.84    & 0.80    & 0.75    & 0.77    & 0.89    & 0.86    & 0.63    & \bf{0.94}    & \bf{0.94}    & \bf{0.87}    \\
   & precision          & 0.86 & \bf{0.80} & 0.46 & 0.72 & 0.63 & 0.49 & \bf{0.88} & 0.70 & \bf{0.98} & 0.86 & 0.71 & 0.94 \\
  & recall           & \bf{0.88} & 0.78 & \bf{0.41} & 0.85 & 0.48 & 0.36 & 0.86 & \bf{0.84} & 0.02 & 0.87 & 0.81 & 0.02 \\
  & F1-score          & \bf{0.86} & \bf{0.79} & \bf{0.39} & 0.78 & 0.54 & \bf{0.39} & \bf{0.86} & 0.74 & 0.01 & \bf{0.86} & 0.73 & 0.01 \\\midrule
\multirow{4}{*}{LNS (Adp vs LCM)} & AUC & \bf{0.79}    & 0.88    & 0.82    & 0.64    & 0.68    & 0.62    & 0.66    & 0.78    & 0.62    & 0.78    & \bf{0.91}    & \bf{0.85}    \\
 & precision & \bf{0.83} & \bf{0.69} & \bf{0.69} & 0.77 & 0.44 & 0.33 & 0.81 & 0.62 & 0.63 & 0.80 & 0.62 & 0.59 \\
 & recall & \bf{0.94} & \bf{0.46} & \bf{0.28} & 0.92 & 0.15 & 0.14 & 0.93 & \bf{0.46} & 0.09 & 0.92 & 0.43 & 0.08 \\
& F1-score & \bf{0.88} & \bf{0.53} & \bf{0.37} & 0.84 & 0.21 & 0.18 & 0.86 & 0.51 & 0.12 & 0.86 & 0.48 & 0.11 \\\midrule
\multirow{4}{*}{LNS (Adp vs 10)}  & AUC & 0.61    & 0.49    & 0.79    & 0.51    & 0.45    & 0.55    & 0.54    & \bf{0.54}    & 0.60    & \bf{0.66}    & 0.43    & \bf{0.89}   \\
& precision & \bf{0.91} & 0.42 & \bf{0.53} & 0.89 & \bf{0.99} & 0.34 & 0.90 & 0.48 & 0.49 & 0.90 & 0.52 & 0.51 \\
& recall  & 0.97 & \bf{0.04} & \bf{0.24} & \bf{0.98} & 0.00 & 0.02 & \bf{0.98} & 0.02 & 0.08 & 0.97 & 0.03 & 0.08 \\
& F1-score  & \bf{0.94} & \bf{0.05} & \bf{0.31} & \bf{0.94} & 0.00 & 0.03 & 0.93 & 0.03 & 0.10 & 0.93 & 0.03 & 0.10 \\\bottomrule 
\end{tabular}
}
\end{table}


\section{Conclusion}
\added[]{Car sequencing problems are NP-hard, and as such the quest for better solution methods is always an ongoing one.}
In this paper, we have performed an instance space analysis for the car sequencing problem, which enabled us to gain deeper insights into problem hardness and algorithm performance. This analysis was achieved by extracting a vector of features to represent an instance and projecting the feature vectors onto a two-dimensional space via principal component analysis. The instance space provided a clear visualization for the distribution of instances, their feature values and the performance of algorithms. We systematically generated instances with controllable properties and showed that the newly generated instances complement the existing benchmark instances well in the instance space. Further, our analysis showed that the instances with a high utilisation or a large average number of options per car class are mostly hard to solve. Finally, we built algorithm selection models using machine learning to select the best algorithm for solving an instance, and found that our new adaptive large neighborhood search algorithm performed the best on hard problem instances, while the MIP solver, Gurobi, performed better on medium-hard instances. Our result based on decision tree provides a simple and effective guideline on which algorithm should be selected to solve a given instance. \added[]{The new benchmark instances created here, and the insight as to which parts of the instance space are most amenable to integer programming versus local search approaches is expected to provide a useful basis for any future research in designing better algorithms for car sequencing problems.}


\section*{Acknowledgement}
This work was partially supported by an ARC Discovery Grant (DP180101170) from Australian Research Council.

\section*{Appendix A: The Pseudocode of LR-ACO}
For the sake of completeness we briefly outline the LR-ACO method in Algorithm~\ref{alg:lr_rcp} here. For details see \citep{Thiruvady2014a}.
The algorithm takes as input, a car sequencing problem instance. The algorithm commences by initializing a best known solution ($\pbs$), the Lagrangian multipliers ($\lambda$) and a number of parameters to be used during the execution of the algorithm. The parameters are (a) scaling factor ($\gamma$), (b) gap to best solution ($gap$), (c) best lower bound ($LB^*$) and best upper bound ($UB^*$) and (d) the pheromone information ($\tau$). Following this, the main part of the algorithm commences, and it runs while three different termination criteria are not met (including the convergence of $\gamma$, a small gap and an iteration limit). Within this loop, there are several function calls, and they are (in order):
\begin{itemize}
\item $x$ = Solve($\lambda^i$, $LB$): solves the relaxed problem, which returns a solution stored in $x$. The solution specifies the positions where each car has been sequenced.
\item $\pi$ = GenerateSequence($x$): $x$ may not necessarily be feasible, and from this infeasible solution a feasible one is generated, which is stored in $\pi$. Note, this is done by selecting a position which has two or more cars and then distributing all excess cars to the nearest empty position (one where no car has been assigned). 
\item ImproveUB($\pi$,$\Tau$): an optional improvement procedure, which takes a sequence and applies ACO to improve it. Note, for this study, we do not use this step as it was previously shown not to be beneficial while using up valuable time. 
\item UpdateBest($\pbs$,$\pi$,$\gamma$): updates the best known solution $\pbs$ if a better one is found. The step size $\gamma$ is reduced if the Lagrangian lower bound does not improve sufficiently.
\item UpdateMult($\lambda^i$, $LB^*$, $UB^*$, $x$, $\gamma$): the Lagrangian multipliers are updated using subgradient optimization as follows: \\$\lambda^i_t := \lambda^{i-1}_t + \gamma \frac{UB^*-LLR(\lambda^i)}{||x||} (\sum_i x_{it} - 1)$.
\end{itemize}

$UB^*$ is set to be the cost of the best known solution  $\pbs$ (Line~11), the final two steps within the loop (Lines~13--14) are to update the gap and iteration count. The last step of the method is to run ACO using the converged pheromone matrix and biasing the search towards the best known solution (Line~16). This was shown to be the most effective approach in obtaining high-quality solutions \citep{Thiruvady2014a}. On completion of the procedure, the best solution is  $\pbs$, which is returned as the output of the algorithm.

\begin{algorithm}[!t]
\caption{LR-ACO}
\small
\label{alg:lr_rcp}
\begin{algorithmic}[1]
\State {\bf \sc input:} Car sequencing instance
\State $\pbs := \NULL$ (best solution)
\State initialize $\lambda^0_{t} = 0, \forall t \in \{1,\ldots,D\}$ 
\State $\gamma := 2.0$, $gap:=\infty$, $UB^* := \infty$, $LB^*:=-\infty$ 
\State $\tau_{ij}=1.0$ $\forall i \in \{1,\ldots,D\}, j \in \{1,\ldots,K\}$
\While{$\gamma > 0.01$ $\&$ $gap>0.01$  $\&$ $i<1000$}
\State $x$ = Solve($\lambda^i$, $LB$)
\State $\pi$ = GenerateSequence($x$)
\State ImproveUB($\pi$,$\Tau$)
\State UpdateBest($\pbs$,$\pi$,$\gamma$)
\State $UB* = f(\pbs)$
\State UpdateMult($\lambda^i$, $LB^*$, $UB^*$, $x$, $\gamma$)
\State $gap = \frac{UB^*-LB^*}{UB^*}$
\State i $\leftarrow$ i + 1
\EndWhile
\State $\pbs \leftarrow$ ACO($\Tau$,$\pbs$)
\State {\bf \sc output:} $\pbs$
\end{algorithmic}
\end{algorithm}

\section*{Appendix B: The Pseudocode of LNS}

\begin{algorithm}[!t]
\caption{{\sf LNS}} \label{alg:ws}
\small
\begin{algorithmic}[1]
\State {\bf \sc input:} $w$, $s$, $t_{\mathrm{lim}}$
\State $\pi :=$ GenerateInitialSequence()

\While{$t_{\mathrm{lim}}$ not exceeded}
\State $\hat{s} \leftarrow$  1
\While{$\hat{s}$ $\leq$ D}
	\State $\hat{\pi} \leftarrow$ Solve($\hat{s}$, $\hat{s}+w$, $\pi$)
	\State $\pi \leftarrow$ Update($\hat{\pi}$)
    \State $\hat{s}$ $\leftarrow$ $\hat{s}$ + $s$
\EndWhile
\EndWhile
\State {\bf \sc output:} $\pi$
\end{algorithmic}
\end{algorithm}

Algorithm~\ref{alg:ws} provides the LNS procedure. The inputs to the algorithm are a window size $w$, a shift size $s$ and a time limit $t_{\mathrm{lim}}$. The first step is to generate a feasible solution (GenerateInitialSequence()), which is stored in $\pi$. There are several possibilities of how this sequence may be found, but as shown in the original study, one Lagrangian subproblem is solved (see Section~\ref{sec:lraco}), the resulting solution is repaired to be feasible, and finally improved by ACO. The next step of the algorithm is to initialize starting point of the sliding window $\hat{s}$. The inner loop executes until the whole sequence has been solved (Line~5). Within this loop, the MIP is solved in the function Solve($\cdot$). Here, the only free variables (associated with the $x$ variables) are those between $\hat{s}$ and $\hat{s}+w$, and $\pi$ is used as a warmstart or seeding solution for the MIP. The output is stored in $\hat{\pi}$, and in the next step (Line~7), the best known solution is updated if a new better one is found. Following this, the starting point of the window is shifted by $s$ positions. The inner loop is repeated for multiple iterations until the time limit is reached. Once the procedure is complete, the algorithm outputs the best solution found $\pi$. 


\bibliographystyle{apalike}
\bibliography{bib2,carsequencing}   

\end{document}